\documentclass{amsart}

\usepackage{amssymb}
\usepackage{amsmath}
\usepackage{graphicx,amsthm,multicol}
\usepackage{hyperref}
\usepackage{epsfig}
\usepackage[hmargin={1.5in,1in},vmargin=1in]{geometry}%,includehead,includefoot]{geometry}
\usepackage{tikz,float,algorithm2e,mathrsfs,nccmath,mathtools, circuitikz, pgfplots}
\usetikzlibrary{positioning}
\usepackage[normalem]{ulem}
\usepackage{makecell}
\usepackage{relsize}
\usepackage{eqparbox}
\usepackage{lineno}
\usepackage{listings}
\usepackage{xcolor} %red, green, blue, yellow, cyan, magenta, black, white
\usepackage{amsaddr}

\RestyleAlgo{ruled}
\SetKwComment{Comment}{/* }{ */}
\usetikzlibrary{decorations.pathmorphing}
\tikzset{
	ragged border/.style={ decoration={random steps, segment length=1mm, amplitude = 0.5mm},
	decorate,
	}
}

\DeclareMathOperator{\sech}{sech}

\definecolor{mygreen}{RGB}{28,172,0} % color values Red, Green, Blue
\definecolor{mylilas}{RGB}{170,55,241}

\makeatletter
\DeclareRobustCommand{\PHP}{%
  \begingroup
  \leavevmode\,\vphantom{P}%
  \dimen\z@=.5\fontcharht\font`P\relax
  \dimen\tw@=0.33333\dimen\z@
  \ooalign{%
    \raisebox{\dimexpr\dimen\z@+2\dimen\tw@-0.4pt}{\rule{\fontcharwd\font`P}{0.4pt}}\cr
    \raisebox{\dimexpr\dimen\z@+\dimen\tw@-0.2pt}{\rule{\fontcharwd\font`P}{0.4pt}}\cr
    P\cr
  }%
  \,\endgroup
}
\makeatother

\usetikzlibrary{arrows,calc,decorations.markings,math,arrows.meta}

\newtheorem{theorem}{Theorem}[section]

\theoremstyle{definition}

\begin{document}

\title[Optimization Models for Dengue Control with Wolbachia]{Balancing Economic Cost and Disease Impact: Optimization Models for Wolbachia-Based Dengue Control} %% Article title

\author{
Kyrho Corum,
Renier Mendoza$^{*}$,\\
Victoria May P. Mendoza,
Arrianne Crystal Velasco
}

\address{
Institute of Mathematics, University of the Philippines Diliman, Quezon City, Philippines
}

\footnotetext{$^{*}$Corresponding author: rmendoza@math.upd.edu.ph}

\maketitle
%\linenumbers

%% Abstract
\begin{abstract}
Dengue, which affects millions of people each year, is one of the most common diseases transmitted by infected \textit{Aedes aegypti} mosquitoes. In the Philippines, the annual economic cost of dengue infections is estimated at around PHP 17 billion. Previous studies have shown that controlling the population of mosquitoes capable of transmitting the dengue virus can effectively reduce dengue infection rates. 
This study explores the use of Wolbachia as a strategy for dengue control by targeting mosquitoes.
Since the release of Wolbachia-infected mosquitoes involves substantial costs, careful planning is necessary to balance disease control with the associated economic burden.
To address this, we propose a mathematical model that captures the dynamics of releasing Wolbachia-carrying mosquitoes and the transmission of dengue in a population. We formulate single- and multi-objective optimization frameworks to minimize the economic costs associated with releasing Wolbachia-infected mosquitoes and the hospitalization costs resulting from dengue infections.  
This study aims to provide insights into the practical application of Wolbachia-based interventions for controlling dengue transmission. While the analysis is grounded in the Philippine context, the approach is general enough to be applicable to other dengue-endemic countries.
\end{abstract}

\section{Introduction}
\label{sec:intro}

Dengue is a viral infection that is passed from mosquitoes to humans \cite{dengue}. As of June 23, 2025, there are 123,291 confirmed dengue cases in the Philippines, with a 0.4\% fatality rate \cite{denguecases}.  On average, the annual economic cost of dengue in the country reaches  PHP 16,947,736,523 \cite{Cj}. This includes reimbursements from the Philippine Health Insurance Corporation, hospitalization costs incurred by patients, and productivity losses due to hospitalization. As a result, various mitigation strategies have been implemented to reduce dengue cases and minimize its economic impact. Usual dengue mitigation strategies include household water container management, fogging, and use of insecticides \cite{denguecontrol}. In 2011, a different dengue suppression strategy was introduced, which makes use of \textit{Aedes} mosquitoes injected with bacteria called Wolbachia \cite{history-aus}.

Wolbachia is a common group of bacteria found in the reproductive tissues of arthropods, transmitted through the cytoplasm of eggs, and has evolved mechanisms to manipulate reproduction of its hosts \cite{wolbachia-intro}. It is known to cause a bendy proboscis, which reduces the biting success of mosquitoes. Also, Wolbachia causes a time delay in the infectiousness of mosquito saliva. This affects viral infections such as dengue, yellow fever, and chikungunya \cite{wolbachia-effects1}. Lastly, Wolbachia causes a reduction in infectious mosquitoes in a population through cytoplasmic incompatibility \cite{ci-factsheet}. Cytoplasmic incompatibility occurs when male mosquitoes with Wolbachia mate with female mosquitoes without Wolbachia. This leads to females laying eggs that would not hatch. On the other hand, if a female mosquito possesses Wolbachia, then regardless of the male's infection to Wolbachia, their offspring will carry Wolbachia \cite{ci-factsheet}. 

The use of Wolbachia to minimize dengue transmission started in 2011, with 298,900 adult mosquitoes released across 374 sites in Cairns, Australia within a span of nine weeks \cite{history-aus}. A follow-up study in 2020 shows a 96\% reduction in dengue incidence \cite{ausresults}. The release of Wolbachia-infected mosquitoes was also tested in Yogakarta, Indonesia in 2014, with a 77\% reduction in dengue incidence \cite{wolbachia-effects1}. In 2016, Singapore launched Project Wolbachia, which initially covered 3,941 households. An 80-90\% reduction in the \textit{Aedes aegypti} mosquito population was reported compared to areas without release, resulting in a 75\% reduction in the likelihood of dengue infection \cite{singaporewolbachia2}. In 2024, the project was expanded to five more areas, covering 480,000 to 580,000 households \cite{singaporewolbachia2,singaporewolbachia}. In the Philippines, interest in using Wolbachia-infected mosquitoes to combat dengue has been expressed since 2020. The Department of Health considered implementing such a project in the Bicol region \cite{wolbachia-ph1}. In 2024, the department is set to conduct research on the use of Wolbachia to reduce dengue transmission in the Philippines \cite{wolbachia-ph2}.

Mathematical models have been used to study the transmission dynamics of dengue. In 1998, the idea of using an SIR-SI model for dengue transmission was introduced, where the human population is divided into susceptible, infected, and recovered populations while the vector population is divided into susceptible and infected populations \cite{SIRSI}. In 2018, de los Reyes V and Escaner introduced an SIJR-SI model for dengue transmission, where the infected human population is divided into healthcare-seeking and nonhealthcare-seeking classes, and applied the model in the Philippine setting \cite{escaner}. In the same year, a deterministic model is proposed where the vector population is divided into immature: comprised of eggs, larvae, and pupae, and adult vector populations. The model also considered vertical transmission, or the transmission of the dengue virus through life stages, and the effects of temperature on the spread of the virus \cite{dengueimmature}.

There are also mathematical models that have studied the spread of Wolbachia infection among mosquitoes. In 2018, Qu et al. looked at the transmission of Wolbachia among the vector population, considering cytoplasmic incompatibility. This involves dividing the vector population according to sex, life stage, and Wolbachia infection status \cite{qumodel}. In 2021, Hu et al. looked at the effects of releasing Wolbachia-infected mosquitoes while considering the effects of pesticides in the general mosquito population \cite{humodel}. Unlike the model of Qu, this model considers only the status of Wolbachia infection and does not distinguish between the vector life stages and sex.

Studies involving Wolbachia in the Philippines focus on the detection of the virus. In 2024, a study suggests the natural prevalence of Wolbachia in Metropolitan Manila \cite{wolbphil}. The study notes that while Wolbachia has been detected in Metro Manila, the prevalence of the virus is low. A study by Sambile shows that one out of 70 samples of \textit{Aedes aegypti} specimens from the provinces of Laguna and Cavite carries Wolbachia.

Optimal control of dengue has also been studied using mathematical models. In 2021, Khan considered an optimal control problem minimizing the number of dengue-infected vectors by employing dengue prevention measures such as the use of mosquito nets and mosquito repellents as well as the use of insecticides against mosquitoes \cite{optifogging}. In 2023, Pongsumpun et al. examined the optimal vaccination against the dengue virus and compared the transmission of dengue under the optimal scheme with the transmission of dengue with the presence of mosquito nets and pesticides \cite{dengueopti}. Moreover, the minimization of the number of dengue-infected individuals by determining the optimal dengue prevention and patient treatment scheme is examined in \cite{dengueopti2}. 
In a related study, de Vasconcelos et al. employed a multi-objective optimization framework to evaluate insecticide-based vector control strategies by minimizing intervention costs and hospitalization expenses \cite{wolbachiaopti}. However, their work does not consider Wolbachia-based mosquito releases, the associated release dynamics and deployment costs. In 2025, Gonzales considered an optimal control problem that minimizes the release of Wolbachia-infected mosquitoes using a mathematical model involving impulsive differential equations to mimic instantaneous release of mosquitoes \cite{wolbachiaopti2}.

In this study, we develop a mathematical model that captures dengue transmission dynamics while incorporating the release of Wolbachia-infected mosquitoes. To the best of our knowledge, no model has incorporated Wolbachia into a dengue model in the Philippines. Different release schemes are considered to study the effect of releasing Wolbachia-infected mosquitoes on the number of dengue-infected individuals. This motivates the formulation of an optimization problem aiming to reduce the costs associated with dengue control and dengue infection. Moreover, this work proposes a multi-objective optimization framework to balance the economic cost of releasing Wolbachia-infected mosquitoes with the goal of minimizing hospitalizations due to dengue, which constitutes the main contribution of this work. While multiobjective optimization frameworks are already used in studies on controlling dengue, there are no published works on the use of such framework to Wolbachia intervention strategies at the time of writing. The findings of this study aim to support policymakers in assessing the feasibility and effectiveness of Wolbachia-based interventions for dengue mitigation.

The rest of the paper is structured as follows. Section \ref{sec:model} discusses the mathematical model that describes the transmission of dengue on the vector and human populations while considering the spread of Wolbachia on the vector population. Section \ref{sec:analysis} presents the existence, uniqueness, positivity, and boundedness of the proposed model. Numerical simulations of different Wolbachia release schemes using the proposed model and formulation of single- and multi-objective optimization problems minimizing the costs associated to dengue are done in Section \ref{sec:numerical}. Section \ref{sec:results} gives the optimal solution to the proposed optimization problems and how they vary after considering the effects of several factors such as production capacity and unit cost. 
Finally, Section \ref{sec:conc} summarizes the main results of the study and highlights some recommendations for future work.

\section{The Proposed Model}
\label{sec:model}
In this section, a mathematical model that describes the transmission of dengue, considering the spread of Wolbachia, is formulated. To do this, we incorporate a Wolbachia population dynamics model to a dengue transmission model. 

We start with the dengue transmission model presented in \cite{escaner}. Here, we look at two clusters: the human cluster and the mosquito cluster. Compartments in the human cluster and in the mosquito cluster are denoted by the subscripts $h$ and $v$, respectively. 

The vector cluster is divided into two: the susceptible compartment ($S_v$) and the infected compartment ($I_v$). The number of susceptible mosquitoes increases at a rate of $b_vN_v\left(1-\dfrac{N_v}{K_a}\right)$, where $b_v$ is the birth rate of mosquitoes, $N_v = S_v+I_v$ is the total population of mosquitoes, and $K_a$ is the carrying capacity of the body of water. Upon biting a human infected with dengue, the mosquitoes move from $S_v$ to $I_v$ compartment at a rate of $BC_{hv}$, where $B$ is the mosquito biting rate and $C_{hv}$ is the dengue transmission rate from human to mosquito. For both the susceptible and infected compartments, we consider a death rate $\mu_v$.

The human cluster is divided into four: the susceptible compartment ($S_h$), the infected compartment ($I_h$), the healthcare-seeking compartment ($J_h$), and the recovered compartment ($R_h$). The number of susceptible humans increase at a rate of $b_hN_h$, where $b_h$ is the birth rate for humans and $N_h = S_h+I_h+J_h+R_h$ is the total population of humans. Being bitten by a female infected mosquito would make susceptible individuals infected at a rate of $BC_{vh}$, with $C_{vh}$ being the transmission rate from vector to human. Only a fraction $\alpha$ of infected individuals seek healthcare upon experiencing symptoms. Healthcare-seeking individuals do not contribute to the transmission of dengue among mosquitoes. After a few days, infected individuals recover at a rate of $\gamma$ for non-healthcare seeking individuals and at a rate of $\theta$ for healthcare-seeking individuals. For all human compartments, we consider a death rate $\mu_h$.

To model the dynamics of dengue transmission in the presence of Wolbachia, we consider dividing the compartments in the vector cluster by the sex of the mosquitoes and by the stage of life of the female mosquitoes. Since Wolbachia-infected mosquitoes have certain characteristics that differ from uninfected mosquitoes, we introduce the following changes. First, we do not divide the male mosquito population into susceptible and infected compartments because only female mosquitoes carry dengue. This is also the case for aquatic-stage mosquitoes. Lastly, since the presence of Wolbachia reduces the transmission of dengue, we consider a transmission rate $C_{vh}^w$ from Wolbachia-infected vectors to humans. Moving forward, the superscript $w$ will be used to denote the Wolbachia-infected population.

The following system of differential equations describes the transmission of dengue with the presence of Wolbachia within the human population:
{\allowdisplaybreaks
\begin{align}
\begin{split}
    \dfrac{dS_h}{dt} &= b_hN_h - \left(BC_{vh}\dfrac{I_v}{N_h}+BC_{vh}^w\dfrac{I_v^w}{N_h}+\mu_h\right)S_h \\
    \dfrac{dI_h}{dt} &= (1-\alpha)\left(BC_{vh}\dfrac{I_v}{N_h}+BC_{vh}^w\dfrac{I_v^w}{N_h}\right)S_h - (\mu_h+\gamma)I_h \\
    \dfrac{dJ_h}{dt} &= \alpha\left(BC_{vh}\dfrac{I_v}{N_h}+BC_{vh}^w\dfrac{I_v^w}{N_h}\right)S_h - (\mu_h+\theta)J_h \\
    \dfrac{dR_h}{dt} &= \gamma I_h+\theta J_h - \mu_h R_h,
\end{split}
\label{modelhuman}
\end{align}}
with $N_h$ being the total human population, $N_v$ being the total vector population, $I_v = I_{vf}+I_{vfp}$ being the total dengue-infected vector population without Wolbachia, and $I_v^w = I_{vf}^w+I_{vfp}^w+I_{vfp}^s$ being the total dengue-infected vector population with Wolbachia. The variable $I$ denotes dengue-infected populations, the subscript $f$ denotes the female vector population, and the subscript $p$ denotes pregnancy of the female vector.

To model the transmission of Wolbachia on the vector population, we adopt the compartmental model presented in \cite{qumodel}, which considers the mating of mosquitoes and the maternal transmission of Wolbachia by dividing the male and female mosquito population into uninfected %($M$, $S_{vf*}$, $I_{vf*}$) 
and Wolbachia-infected.  
%($M^w$, $S_{vf*}^w$, $I_{vf*}^w$).

Uninfected females ($S_{vf}$, $I_{vf}$) and Wolbachia-infected females ($S_{vf}^w$, $I_{vf}^w$) mate with Wolbachia-infected males ($M^w$) or with uninfected males ($M$) and enter the pregnancy stage at a mating rate $\sigma$. We define \[m(M,M^w) = \dfrac{M}{M+M^w} \quad \text{and} \quad m_w = 1-m\] to be the probability that a male is uninfected, and the probability that a male is infected, respectively. Pregnant females can either be uninfected ($S_{vfp}$, $I_{vfp}$) when an uninfected male mates with an uninfected female at a rate of $\sigma\cdot m$, Wolbachia-infected ($S_{vfp}^w$, $I_{vfp}^w$) when an infected female mates with a male at a rate of $\sigma$, or sterile ($S_{vfp}^s$, $I_{vfp}^s$) when an uninfected female mates with an infected male at a rate of $\sigma\cdot m_w$. The rates $\sigma\cdot m$ and $\sigma\cdot m_w$ are consistent with the compatibility of the transmission of Wolbachia as described in \cite{ci-factsheet}.

Pregnant females will eventually produce aquatic-stage mosquitoes. Uninfected pregnant females will only produce uninfected eggs ($A$) at a rate of $\eta$, sterile pregnant females do not produce eggs, and Wolbachia-infected pregnant females produce infected offspring ($A^w$) at a fraction of $v_w$ and a rate of $\eta_w$, and uninfected offspring at a fraction of $v$ and a rate of $\eta$. In particular, we define the oviposition rates, or the average number of eggs produced during a female mosquito's lifespan, for Wolbachia-uninfected and Wolbachia-infected females as \[\eta(A,A^w)=\phi\left(1-\dfrac{A+A^w}{K_a}\right) \quad \text{and} \quad \eta_w(A,A^w)=\phi_w\left(1-\dfrac{A+A^w}{K_a}\right),\] respectively, where $\phi$ and $\phi_w$ are the egg-laying rates for uninfected and infected mosquitoes, respectively, and $K_a$ is the carrying capacity of the body of water under consideration. The aquatic-stage mosquitoes eventually hatch and grow into adult forms at a rate of $\psi$, where a fraction $b_m$ are males and a fraction $b_f$ are females.

To account for mosquito mortality, we consider the following death rates: $\mu_a$ for aquatic-stage mosquitoes, $\mu_{fu}$ for female uninfected mosquitoes, $\mu_{fw}$ for female infected mosquitoes, $\mu_{mu}$ for male uninfected mosquitoes, and $\mu_{mw}$ for male Wolbachia-infected mosquitoes. Finally, to introduce Wolbachia to the vector population, we consider a release $r(t)$ of Wolbachia-infected aquatic-stage mosquitoes. We have the following equations for the male and aquatic-stage mosquito population:
\begin{align}
    \begin{split}
    \dfrac{dM_v^w}{dt} &= \psi b_m A^w - \mu_{m}^wM_v^w \\
    \dfrac{dM_v}{dt} &= \psi b_m A - \mu_{m} M_v \\
    \dfrac{dA^w}{dt} &= r(t)+ \eta^w v^w (S_{vfp}^w+I_{vfp}^w) - (\psi+\mu_a)A^w \\
    \dfrac{dA}{dt} &= \eta(S_{vfp}+I_{vfp})+\eta^w v(S_{vfp}^w+I_{vfp}^w)-(\psi+\mu_a)A.
    \end{split}
    \label{modelvec1}
\end{align}
Furthermore, the remaining equations describe the vector population in the susceptible and infected compartments with respect to both dengue and Wolbachia. Note that $S$ and $I$ denote susceptibility and infection to dengue, respectively.
{\allowdisplaybreaks
\begin{align}
   \begin{split}
    \dfrac{dS_{vf}^w}{dt} &= \psi b_f A^w - \left(BC_{hv}\dfrac{I_h}{N_h}+\sigma+\mu_{f}^w\right)S_{vf}^w \\
    \dfrac{dS_{vf}}{dt} &= \psi b_f A - \left(BC_{hv}\dfrac{I_h}{N_h}+\sigma+\mu_{f}\right)S_{vf} \\
    \dfrac{dI_{vf}^w}{dt} &= BC_{hv}\dfrac{I_h}{N_h}S_{vf}^w - (\sigma+\mu_{f}^w)I_{vf}^w \\
    \dfrac{dI_{vf}}{dt} &= BC_{hv}\dfrac{I_h}{N_h}S_{vf} - (\sigma+\mu_{f})I_{vf} \\
    \dfrac{dS_{vfp}^w}{dt} &= \sigma S_{vf}^w - \left(BC_{hv}\dfrac{I_h}{N_h}+\mu_{f}^w\right)S_{vfp}^w \\
    \dfrac{dS_{vfp}}{dt} &= \sigma m S_{vf} - \left(BC_{hv}\dfrac{I_h}{N_h}+\mu_{f}\right)S_{vfp} \\
    \dfrac{dS_{vfp}^s}{dt} &= \sigma m^w S_{vf} - \left(BC_{hv}\dfrac{I_h}{N_h}+\mu_{f}\right)S_{vfp}^s \\
    \dfrac{dI_{vfp}^w}{dt} &= \sigma I_{vf}^w + BC_{hv}\dfrac{I_h}{N_h} S_{vfp}^w - \mu_{f}^w I_{vfp}^w \\
    \dfrac{dI_{vfp}}{dt} &= \sigma m I_{vf} + BC_{hv}\dfrac{I_h}{N_h} S_{vfp} - \mu_{f} I_{vfp} \\
    \dfrac{dI_{vfp}^s}{dt} &= \sigma m^w I_{vf} + BC_{hv}\dfrac{I_h}{N_h}S_{vfp}^s - \mu_{f}^w I_{vfp}^s.
    \end{split}
    \label{modelvec2}
\end{align}
}
We note an important assumption that the right-hand side of the differential equations in systems \eqref{modelhuman}-\eqref{modelvec2} are continuous functions on $[0,T]$ for some $T\in\mathbb{R}$.

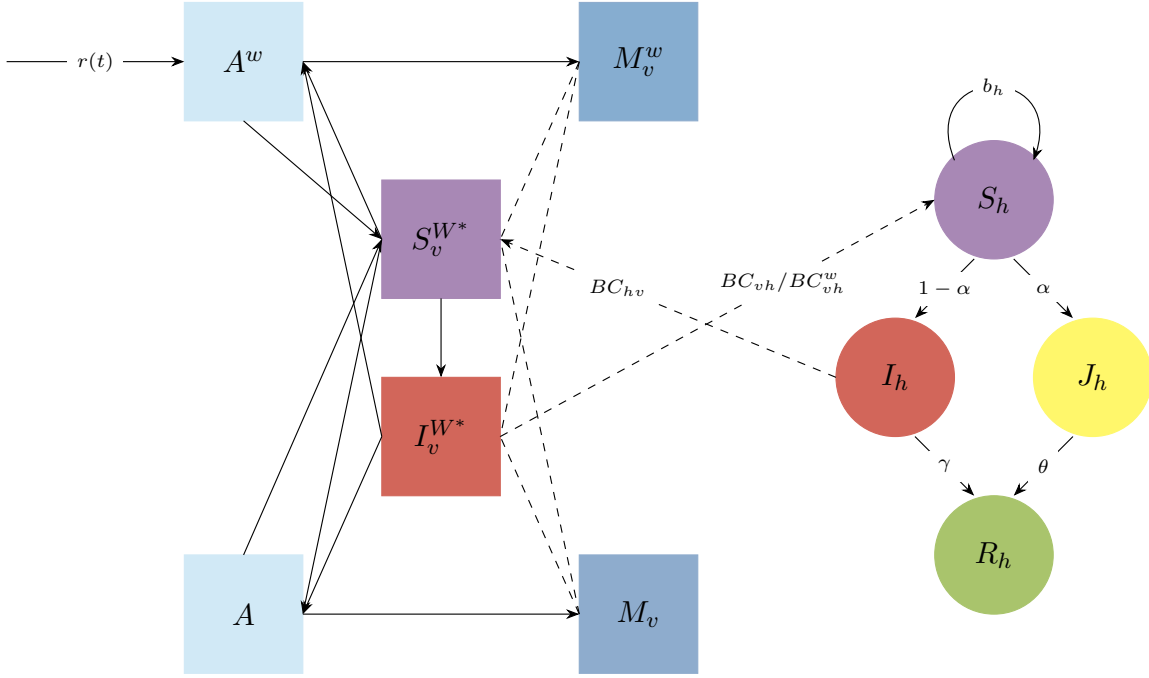
\begin{figure}
\centering
\resizebox{\textwidth}{!}{%
\begin{circuitikz}
\tikzstyle{every node}=[font=\large]
\draw [ color={rgb,255:red,168; green,136; blue,181} , fill={rgb,255:red,168; green,136; blue,181}] (19.5,12) circle (0.75cm);
\node [font=\large] at (19.5,12) {$S_h$};
\draw [ color={rgb,255:red,210; green,102; blue,90} , fill={rgb,255:red,210; green,102; blue,90}] (18.25,9.75) circle (0.75cm);
\node [font=\large] at (18.25,9.75) {$I_h$};
\draw [ color={rgb,255:red,254; green,247; blue,108} , fill={rgb,255:red,255; green,247; blue,107}] (20.75,9.75) circle (0.75cm);
\node [font=\large] at (20.75,9.75) {$J_h$};
\draw [ color={rgb,255:red,169; green,196; blue,108} , fill={rgb,255:red,169; green,196; blue,108}] (19.5,7.5) circle (0.75cm);
\node [font=\large] at (19.5,7.5) {$R_h$};
\draw [->, >=Stealth] (19.75,11.25) -- (20.5,10.5)node[pos=0.5, fill=white]{\scriptsize{$\alpha$}};
\draw [->, >=Stealth] (19.25,11.25) -- (18.5,10.5)node[pos=0.5, fill=white]{\scriptsize{$1-\alpha$}};
\draw [->, >=Stealth] (20.5,9) -- (19.75,8.25)node[pos=0.5, fill=white]{\scriptsize{$\theta$}};
\draw [->, >=Stealth] (18.5,9) -- (19.25,8.25)node[pos=0.5, fill=white]{\scriptsize{$\gamma$}};
\draw [->, >=Stealth] (19,12.5) .. controls (18.5,13.75) and (20.5,13.75) .. (20,12.5) node[pos=0.5, fill=white]{\scriptsize{$b_h$}};
\draw [ color={rgb,255:red,168; green,136; blue,181} , fill={rgb,255:red,168; green,136; blue,181}] (13.25,12.25) rectangle (11.75,10.75);
\node [font=\large] at (12.5,11.5) {$S_v^{W^*}$};
\draw [ color={rgb,255:red,210; green,102; blue,90} , fill={rgb,255:red,210; green,102; blue,90}] (13.25,9.75) rectangle (11.75,8.25);
\node [font=\large] at (12.5,9) {$I_v^{W^*}$};
\draw [ color={rgb,255:red,142; green,172; blue,205} , fill={rgb,255:red,134; green,173; blue,208}] (15.75,14.5) rectangle (14.25,13);
\node [font=\large] at (15,13.75) {$M_v^w$};
\draw [ color={rgb,255:red,209; green,233; blue,246} , fill={rgb,255:red,209; green,233; blue,246}] (10.75,14.5) rectangle (9.25,13);
\node [font=\large] at (10,13.75) {$A^w$};
\draw [ color={rgb,255:red,209; green,233; blue,246} , fill={rgb,255:red,209; green,233; blue,246}] (10.75,7.5) rectangle (9.25,6);
\node [font=\large] at (10,6.75) {$A$};
\draw [ color={rgb,255:red,142; green,172; blue,205} , fill={rgb,255:red,142; green,172; blue,205}] (15.75,7.5) rectangle (14.25,6);
\node [font=\large] at (15,6.75) {$M_v$};
\draw [->, >=Stealth, dashed] (17.5,9.75) -- (13.25,11.5)node[pos=0.65, fill=white]{\scriptsize{$BC_{hv}$}};
\draw [->, >=Stealth, dashed] (13.25,9) -- (18.75,12)node[pos=0.65, fill=white]{\scriptsize{$BC_{vh}/BC_{vh}^w$}};
\draw [dashed] (14.25,13.75) -- (13.25,11.5);
\draw [dashed] (14.25,13.75) -- (13.25,9);
\draw [dashed] (14.25,6.75) -- (13.25,9);
\draw [dashed] (14.25,6.75) -- (13.25,11.5);
\draw [->, >=Stealth] (11.75,11.5) -- (10.75,13.75);
\draw [->, >=Stealth] (11.75,11.5) -- (10.75,6.75);
\draw [->, >=Stealth] (11.75,9) -- (10.75,6.75);
\draw [->, >=Stealth] (11.75,9) -- (10.75,13.75);
\draw [->, >=Stealth] (10,13) -- (11.75,11.5);
\draw [->, >=Stealth] (10,7.5) -- (11.75,11.5);
\draw [->, >=Stealth] (10.75,13.75) -- (14.25,13.75);
\draw [->, >=Stealth] (10.75,6.75) -- (14.25,6.75);
\draw [->, >=Stealth] (7,13.75) -- (9.25,13.75)node[pos=0.5, fill=white]{\scriptsize{$r(t)$}};
\draw [->, >=Stealth] (12.5,10.75) -- (12.5,9.75);
\end{circuitikz}}
\caption{Schematic diagram of the dengue transmission model with the presence of Wolbachia. Humans who are susceptible to dengue ($S_h$) become infected ($I_h,J_h$) upon contact with mosquitoes that carry dengue ($I_v^W$) at a rate of $BC_{vh}$ or $BC_{vh}^w$. Mosquitoes who are susceptible to dengue ($S_v^W$) become infected ($I_v^W$) upon contact with nonhealthcare-seeking humans. Those in $I_h$ recover at a rate of $\gamma$ while those in $J_h$ recover at a rate of $\theta$. \\ $^*$Transmission dynamics for Wolbachia in $S_v^W$ and $I_v^W$ is demonstrated in Figure \ref{transmission}.}
\label{finalmodel}
\end{figure}

\begin{figure}
\centering
\resizebox{1\textwidth}{!}{%
\begin{circuitikz}
\tikzstyle{every node}=[font=\scriptsize]

\draw [ color={rgb,255:red,142; green,172; blue,205} , fill={rgb,255:red,134; green,173; blue,208}] (22,16) rectangle (20.5,14.5);
\node [font=\large] at (21.25,15.25) {$M_v^w$};
\draw [ color={rgb,255:red,209; green,233; blue,246} , fill={rgb,255:red,209; green,233; blue,246}] (10.75,16) rectangle (9.25,14.5);
\node [font=\large] at (10,15.25) {$A^w$};
\draw [ color={rgb,255:red,209; green,233; blue,246} , fill={rgb,255:red,209; green,233; blue,246}] (10.75,4.75) rectangle (9.25,3.25);
\node [font=\large] at (10,4) {$A$};
\draw [ color={rgb,255:red,142; green,172; blue,205} , fill={rgb,255:red,142; green,172; blue,205}] (22,4.75) rectangle (20.5,3.25);
\node [font=\large] at (21.25,4) {$M_v$};
\draw [->, >=Stealth] (10.75,15.25) -- (20.5,15.25)node[pos=0.5, fill=white]{$b_m\psi$};
\draw [->, >=Stealth] (10.75,4) -- (20.5,4)node[pos=0.5, fill=white]{$b_m\psi$};
\draw [->, >=Stealth] (7,15.25) -- (9.25,15.25)node[pos=0.5, fill=white]{$r(t)$};
\draw [ color={rgb,255:red,212; green,180; blue,208} , fill={rgb,255:red,212; green,180; blue,208}] (18,12.25) rectangle (16.5,10.75);
\draw [ color={rgb,255:red,212; green,180; blue,208} , fill={rgb,255:red,212; green,180; blue,208}] (18,9.75) rectangle (16.5,8.25);
\node [font=\large] at (17.25,11.5) {$S_{vf}^w$};
\node [font=\large] at (17.25,9) {$S_{vf}$};
\draw [ color={rgb,255:red,212; green,137; blue,208} , fill={rgb,255:red,212; green,137; blue,208}] (14.5,13.25) rectangle (13,11.75);
\node [font=\large] at (13.75,12.5) {$S_{vfp}^w$};
\draw [ color={rgb,255:red,212; green,137; blue,208} , fill={rgb,255:red,212; green,137; blue,208}] (14.5,11) rectangle (13,9.5);
\node [font=\large] at (13.75,10.25) {$S_{vfp}^s$};
\draw [ color={rgb,255:red,212; green,137; blue,208} , fill={rgb,255:red,212; green,137; blue,208}] (14.5,8.75) rectangle (13,7.25);
\node [font=\large] at (13.75,8) {$S_{vfp}$};
\draw [dashed] (21.25,14.5) -- (18,11.5);
\draw [dashed] (21.25,14.5) -- (18,9);
\draw [dashed] (21.25,4.75) -- (18,9);
\draw [dashed] (21.25,4.75) -- (18,11.5);
\draw [->, >=Stealth] (16.5,9) -- (14.5,10.25)node[pos=0.5, fill=white]{$\sigma m_w$};
\draw [->, >=Stealth] (16.5,9) -- (14.5,8)node[pos=0.5, fill=white]{$\sigma m$};
\draw [->, >=Stealth] (16.5,11.5) -- (14.5,12.5)node[pos=0.5, fill=white]{$\sigma$};
\draw [->, >=Stealth] (13,8) -- (10.75,4)node[pos=0.5, fill=white]{$\eta$};
\draw [->, >=Stealth] (13,12.5) -- (10.75,4)node[pos=0.5, fill=white]{$\eta_w v$};
\draw [->, >=Stealth] (13,12.5) -- (10.75,15.25)node[pos=0.5, fill=white]{$\eta_w v_w$};
\draw [->, >=Stealth] (10.75,15.25) -- (17.25,12.25)node[pos=0.5, fill=white]{$b_f\psi$};
\draw [->, >=Stealth] (10.75,4) -- (17.25,8.25)node[pos=0.5, fill=white]{$b_f\psi$};
\end{circuitikz}
}%

\label{transmission}
\caption{Transmission dynamics of Wolbachia on $S_v^W$. We note that similar dynamics are observed on $I_v^W$. Mosquitoes in the $S_{vf}^w$ compartment transition into $S_{vfp}^w$ at a rate of $\sigma$. Mosquitoes in the $S_{vf}$ compartment transition into either $S_{vfp}$ or $S_{vfp}^s$ at a rate of $\sigma m$ or $\sigma m_w$, respectively. Mosquitoes in the $S_{vfp}$ compartment lay uninfected eggs ($A$) at a rate of $\eta$ while mosquitoes in the $S_{vfp}^w$ lay either Wolbachia-infected ($A^w$) or uninfected eggs at a rate of $\eta_w v_w$ or $\eta_w v$, respectively. Wolbachia-infected eggs eventually develop into Wolbachia-infected male ($M_v^w$) and (dengue susceptible) infected nonpregnant female mosquitoes while Wolbachia-uninfected eggs develop into Wolbachia-uninfected male ($M_v$) and (dengue susceptible) nonpregnant female mosquitoes. We consider a release $r(t)$ of Wolbachia-infected eggs.}
\end{figure}
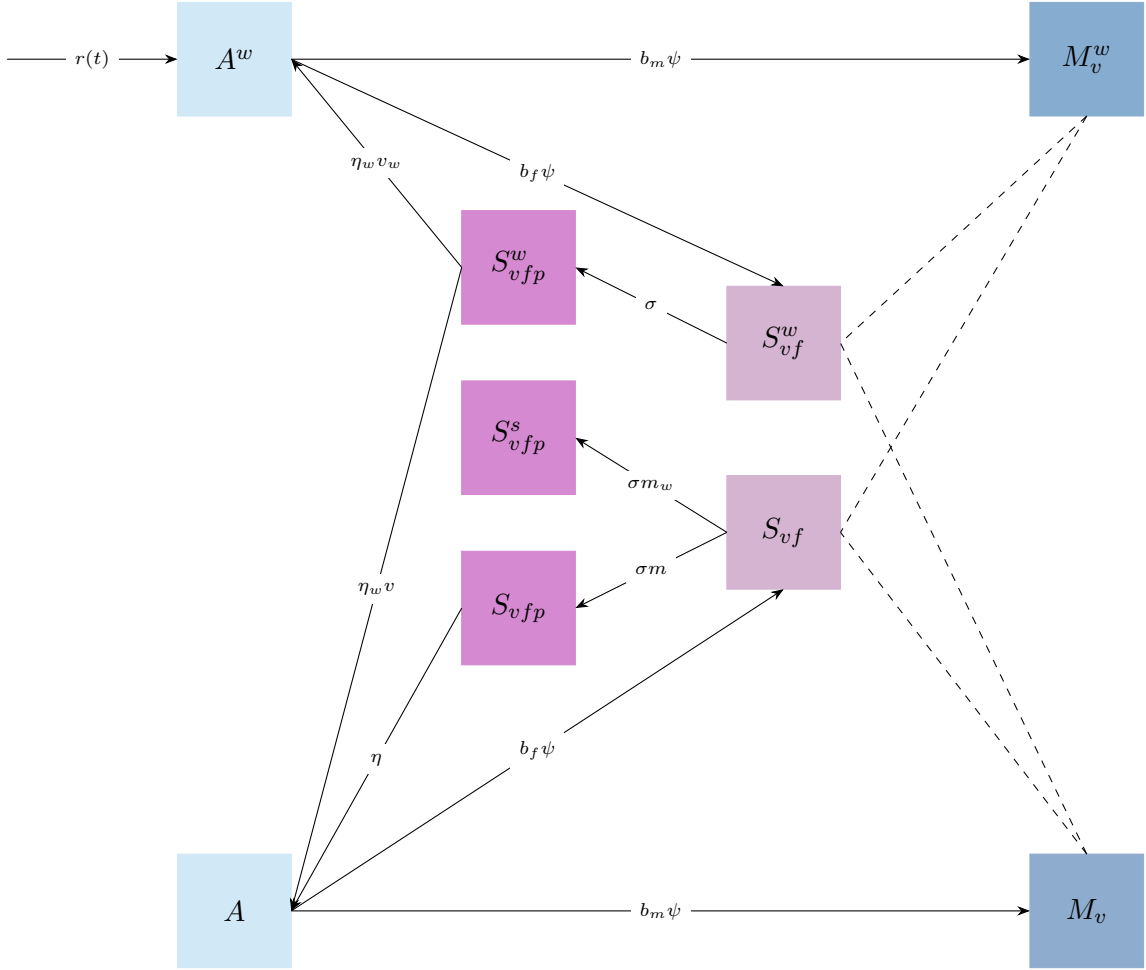

Figure \ref{finalmodel} summarizes the transmission of dengue from mosquito to human, and vice versa, with the presence of Wolbachia. The transmission dynamics of Wolbachia in the $S_v^W$ and $I_v^W$ compartments are detailed in Figure \ref{transmission}.

%\newpage
%\newpage

\newpage

\section{Analysis of the model}
\label{sec:analysis}
In this section, the existence, uniqueness, and positivity or forward invariance of the solution of the proposed model \eqref{modelhuman}-\eqref{modelvec2} is established. 

%In particular, we prove the existence and uniqueness of the solution to the system in Theorem \ref{existence} for $t\geq 0$.
\begin{theorem}[Existence and uniqueness of solutions]
\label{existence}
    There exists a {local} unique solution to the equations in \eqref{modelhuman}-\eqref{modelvec2} {on the interval $[0,T]$ for some $T\in \mathbb{R}$}. 
    \begin{proof}
         {Let $I$ be a neighborhood for $t_0\in [0,T]$. Consider the mapping $x(t): [0,T] \to \mathbb{R}^{18}$ such that} 
        \begin{align*}
            x(t) = & \left [ S_h(t), I_h(t), J_h(t), R_h(t), M_v^w(t), M_v(t), S_{vf}^w(t), S_{vf}(t), S_{vfp}^w(t),  \right . \\
        & \left . S_{vfp}(t), S_{vfp}^s(t), I_{vf}^w(t), I_{vf}(t), I_{vfp}(t), I_{vfp}^s(t), I_{vfp}^w(t), A^w(t), A(t) \right ]^T.
        \end{align*} 
        We write the initial value problem from System \eqref{modelhuman}-\eqref{modelvec2} as 
        \[\dot{x} = F(t,x(t)), \quad x(0) = x_0,\] 
        %where $x(t)= [ S_h(t), I_h(t), J_h(t), R_h(t), M_v^w(t), M_v(t), S_{vf}^w(t), S_{vf}(t), S_{vfp}^w(t), S_{vfp}(t),$ \\ $ S_{vfp}^s(t), I_{vf}^w(t), I_{vf}(t), I_{vfp}(t), I_{vfp}^s(t), I_{vfp}^w(t), A^w(t), A(t) ]^T \in \mathbb{R}^{18}$ defines a mapping from $[0,+\infty)$ to $\mathbb{R}^{18}$. 
        where $F(t,x(t)) = [ F_1(t,x(t)),\ldots,F_{18}(t,x(t) ]^T$ is a vector-valued function 
        %from $\mathbb{R}^{18}$ to $\mathbb{R}^{18}$ 
        corresponding to the right-hand side of the equations in \eqref{modelhuman}-\eqref{modelvec2}. From its form, we deduce that $F$ is continuous. Furthermore, the entries of the Jacobian matrix of $F$, as shown in Appendix \ref{JacobianF}, are all continuous. Hence, $F$ has continuous first partial derivatives with respect to $x$, 
        %in $\mathbb{R}^{18}$, 
        implying that $F$ is locally Lipschitz in $x$. Thus, by Picard-Lindel\"of theorem \cite{existence}, System \eqref{modelhuman}-\eqref{modelvec2} has a unique solution for {all $t\in I$}. 
    \end{proof}
\end{theorem}

Next, we establish the \textit{forward invariance} of a domain $\mathscr{D} \subset \mathbb{R}^{18}$. This assures that, along with Theorem \ref{existence}, the solution to the system of equations \eqref{modelhuman}-\eqref{modelvec2} remains in $\mathscr{D}$ {for all $t\in I$, where $I$ is a neighborhood of $t_0\in [0,T]$.} 

\begin{theorem}[Forward invariance]
    Let $x = (S_h, I_h, J_h, R_h, M_v^w, M_v, S_{vf}^w, S_{vf}, S_{vfp}^w, S_{vfp},$ \\ $ S_{vfp}^s, I_{vf}^w, I_{vf}, I_{vfp}, I_{vfp}^s, I_{vfp}^w, A^w, A)^T$.
    Assuming that the initial condition $x_0=x(0)$ lies in the epidemiological domain \[\mathscr{D} = \left\{x\in\mathbb{R}_+^{18} \ : x\text{ satisfies }\eqref{upperbound} \right\},\] with the bounds for the vector population being as follows for any $t$:
    \begin{align}
    \begin{split}
        0&\leq A+A^w\leq K_a \\
        0&\leq M_v+M_v^w\leq b_m\dfrac{\psi}{\mu_{mu}}K_a \\
        0&\leq S_{vf}+S_{vf}^w+I_{vf}+I_{vf}^w \leq b_f \dfrac{\psi}{\sigma+\mu_{fu}}K_a \\
        0&\leq S_{vfp}+S_{vfp}^w+S_{vfp}^s+I_{vfp}+I_{vfp}^w+I_{vfp}^s\leq b_f\dfrac{\sigma}{\sigma+\mu_{fu}}\dfrac{\psi}{\mu_{fu}}K_a,
    \end{split}
    \label{upperbound}
\end{align}
the system of equations \eqref{modelhuman}-\eqref{modelvec2} has a unique solution that remains in $\mathscr{D}$ for all $t\in I$, where $I$ is a neighborhood of $t_0\in [0,T]$.

    \begin{proof}
        Let $t_0\in [0,T]$ and consider a neighborhood $I$ of $t_0$. To establish that the solution of the system \eqref{modelhuman}-\eqref{modelvec2} is nonnegative, we argue that for each equation of the form $\dfrac{dy_i}{dt}=f(y),$ for $i=1,\ldots,n$, where $y=(y_1,\ldots,y_n) \in \mathscr{D}$, if $y_i=0$, then $f(y)\geq 0$. Indeed, from Equation \eqref{modelhuman}
        \[S_h=0 \Rightarrow \dfrac{dS_h}{dt} = b_h(I_h+J_h+R_h) \geq 0\] and \[I_h = 0 \Rightarrow \dfrac{dI_h}{dt} = (1-\alpha)\left(BC_{vh}\dfrac{I_v}{S_h+J_h+R_h}+BC_{vh}^w\dfrac{I_v^w}{S_h+J_h+R_h}\right)S_h \\ \geq 0.\] Similar arguments are used for $J_h$ and $R_h$. 
        
        Next, we show the nonnegativity of $M_v^w$, with $M_v$ following the same argument.
        \[M_v^w = 0 \Rightarrow \dfrac{dM_v^w}{dt} = \psi b_m A^w \geq 0.\] 
        To show the nonnegativity of $S_{vf}^w$ and $I_{vf}^w$, we have
        \[S_{vf}^w = 0 \Rightarrow \dfrac{dS_{vf}^w}{dt} = \psi b_f A^w \geq 0\] 
        and 
        \[I_{vf}^w = 0 \Rightarrow \dfrac{dI_{vf}^w}{dt} = BC_{hv} \dfrac{I_h}{N_h}S_{vf}^w \geq 0.\] 
        Again, we use similar arguments to show the nonnegativity of $S_{vf}$ and $I_{vf}$. Now, for the nonnegativity of $S_{vfp}^w$ and $I_{vfp}^w$, we get
        \[S_{vfp}^w = 0 \Rightarrow \dfrac{dS_{vfp}^w}{dt} = \sigma S_{vf}^w \geq 0\] and \[I_{vfp}^w = 0 \Rightarrow \dfrac{dI_{vfp}^w}{dt} = \sigma I_{vf}^w + BC_{hv}\dfrac{I_h}{N_h}S_{vfp}^w \geq 0.\] Proving the nonnegativity of $S_{vfp}$, $S_{vfp}^s$, $I_{vfp}$, and $I_{vfp}^s$ follows the same reasoning. Lastly, \[A^w = 0 \Rightarrow \eta_w v_w (S_{vfp}^w+I_{vfp}^w)+r(t) \geq 0,\] with nonnegativity of $A$ following the same justification.
        
        Furthermore, we show that if the upper bounds of the mosquito population are met, then the mosquito population will not increase further. For example, we show that if $A+A^w=K_a$, then $\dfrac{dA}{dt}+\dfrac{dA^w}{dt}\leq0$. Note that the death rate for Wolbachia-infected mosquitoes is higher than those of uninfected mosquitoes \cite{qumodel}. This is important in establishing the desired inequalities.

        We begin with the death rate for the aquatic stage mosquitoes.
        \[ A+A^w=K_a \Rightarrow \dfrac{dA}{dt}+\dfrac{dA^w}{dt} = r(t) -(\psi+\mu_a)(A+A^w) \leq r(t) -(\psi+\mu_a)K_a.\]  Imposing the restriction $r(t)\leq (\psi+\mu_a)K_a$ assures that $\dfrac{dA}{dt}+\dfrac{dA^w}{dt} \leq 0$. This suggests that the amount of release depends on the carrying capacity for mosquitoes.
        
        For the male mosquitoes, we get $\mu_m \leq \mu_m^w$ that $M_v+M_v^w=\dfrac{b_m\psi}{\mu_{m}}K_a$ implies
        \begin{align*}
\dfrac{dM_v}{dt}+\dfrac{dM_v^w}{dt} &= \psi b_m(A+A^w)-\mu_{m}M_v-\mu_{m}^wM_v^w \\ &\leq \psi b_m(A+A^w) - \mu_{m}(M_v+M_v^w),
        \end{align*}
        with the inequality following from the fact that the mortality rate for Wolbachia-infected mosquitoes are higher than that of the uninfected mosquitoes \cite{qumodel}.
        Using the assumption that $\mu_m(M_v+M_v^w)=\psi b_mK_a$ and the fact that $A+A_w\leq K_a$, we see that 
        \[\dfrac{dM_v}{dt}+\dfrac{dM_v^w}{dt}\leq \psi b_m (A+A^w-K_a) \leq 0.\]
        For the nonpregnant female mosquitoes, see that if $S_{vf}+S_{vf}^w+I_{vf}+I_{vf}^w = \dfrac{b_f\psi}{\sigma+\mu_{f}} K_a$, then
        \begin{align*}
            \dfrac{dS_{vf}}{dt}+\dfrac{dS_{vf}^w}{dt}+\dfrac{dI_{vf}}{dt}+\dfrac{dI_{vf}^w}{dt} &= \psi b_f(A+A^w)-(\sigma+\mu_{fw})(S_{vf}^w+I_{vf}^w) -(\sigma+\mu_{f})(S_{vf}+I_{vf}) \\ &\leq \psi b_f(A+A^w) - (\sigma+\mu_{f})(S_{vf}+I_{vf}+S_{vf}^w+I_{vf}^w). \\ &= \psi b_f (A+A^w-K_a) \leq 0.
        \end{align*}
        Note that the inequality follows from the assumption that $\mu_f\leq \mu_f^w$ \cite{qumodel}. Now, using the assumption that $(\sigma+\mu_f)(S_{vf}+S_{vf}^w+I_{vf}+I_{vf}^w)=\psi b_fK_a$, together with the fact that $A+A_w\leq K_a$, we see that 
        \[\dfrac{dS_{vf}}{dt}+\dfrac{dS_{vf}^w}{dt}+\dfrac{dI_{vf}}{dt}+\dfrac{dI_{vf}^w}{dt} \leq \psi b_f (A+A^w-K_a) \leq 0.\]
        {Lastly, for the pregnant mosquitoes, if $S_{vfp}+I_{vfp}+S_{vfp}^w+I_{vfp}^w+S_{vfp}^s+I_{vfp}^s = \dfrac{\sigma}{\sigma+\mu_{f}}\dfrac{b_f\psi}{\mu_{f}}K_a$, then}
        \begin{align*}
              &\dfrac{dS_{vfp}}{dt}+\dfrac{dI_{vfp}}{dt}+\dfrac{dS_{vfp}^w}{dt}+\dfrac{dI_{vfp}^w}{dt}+\dfrac{dS_{vfp}^s}{dt}+\dfrac{dI_{vfp}^s}{dt} \\ &= \sigma(S_{vf}^w+I_{vf}^w) + \sigma m(S_{vf}+I_{vf}) - \mu_{f}^w(S_{vfp}^w+I_{vfp}^w) - \mu_{f}(S_{vfp}+I_{vfp}+S_{vfp}^s+I_{vfp}^s) \\ &\leq \sigma(S_{vf}^w+S_{vf}+I_{vf}^w+I_{vf}) -\mu_{f}(S_{vfp}+I_{vfp}+S_{vfp}^w+I_{vfp}^w+S_{vfp}^s+I_{vfp}^s). \\ &\leq \sigma \dfrac{b_f\psi}{\sigma+\mu_{f}}K_a - \mu_{f}\dfrac{\sigma}{\sigma+\mu_{f}}\dfrac{b_f\psi}{\mu_{f}}K_a = 0.
            \end{align*}
            {The inequality follows from the fact that $m\leq 1$ and $\mu_f\leq \mu_f^w$. Now, since $(\sigma+\mu_f)(S_{vf}+I_{vf}+S_{vf}^w+I_{vf}^w) \leq b_f \psi K_a$ and $\mu_{f}(S_{vfp}+I_{vfp}+S_{vfp}^w+I_{vfp}^w+S_{vfp}^s+I_{vfp}^s) = \sigma \dfrac{b_f\psi}{\sigma+\mu_f}K_a$, it follows that} \[\dfrac{dS_{vfp}}{dt}+\dfrac{dI_{vfp}}{dt}+\dfrac{dS_{vfp}^w}{dt}+\dfrac{dI_{vfp}^w}{dt}+\dfrac{dS_{vfp}^s}{dt}+\dfrac{dI_{vfp}^s}{dt} \leq \sigma \dfrac{b_f\psi}{\sigma+\mu_{f}}K_a - \sigma \dfrac{b_f\psi}{\sigma+\mu_f}K_a=0.\]
            {These show that if the upper bounds prescribed in \eqref{upperbound} are met, then the mosquito population does not increase. Hence, if the initial condition lies in the domain $\mathscr{D}$, then the solution remains in $\mathscr{D}$ for all $t\in I$.}
    \end{proof}
\end{theorem}

\section{Numerical Methods}
\label{sec:numerical}
In the first subsection, we examine different release schemes $r(t)$ and their impact on the number of dengue infections. The model given by the system of differential equations \eqref{modelhuman}-\eqref{modelvec2} is solved numerically.  In the following two subsections, we present the single- and multi-objective optimization problems. The optimization problems aim to determine $r(t)$ that minimize the cost of the Wolbachia intervention and the societal cost of dengue-related hospitalizations.

\subsection{Simulation of the model} \label{simulation}
To simulate the model, we solve the system of equations in \eqref{modelhuman}- \eqref{modelvec2} using the Matlab built-in function \texttt{ode45}, which is based on a Runge-Kutta method with adaptive step size \cite{odesolver}. We note here that all simulations presented in this work were run on MATLAB R2024b using a 16GB M3 Macbook Air. The parameter values and the initial conditions used in the simulation are listed in Tables \ref{compartmentnames} and  \ref{parameternames}.

\begin{table}
    \centering
    \begin{tabular}{|r|l|c|c|}
        \hline
        Compartment & Description & Initial Value \\
        \hline
        \hline
        $S_h$ & Susceptible humans & 50~000~000 \\
        $I_h$ & Infected humans & 15~000  \\
        $J_h$ & Healthcare-seeking humans & 1~500\\
        $R_h$ & Recovered humans & 5~000  \\
        \hline
        \hline
        $M_v$ & Male vector without Wolbachia & 10~000~000 \\
        $M_v^w$ & Male vector with Wolbachia & 0 \\
        \hline
        $S_{vf}$ & Susceptible female vector &  \\ &without Wolbachia & 7~500~000 \\
        $S_{vf}^w$ & Susceptible female vector &  \\ & with Wolbachia & 0\\
        $S_{vfp}$ & Susceptible pregnant female vector  & \\ & without Wolbachia & 2~500~000  \\
        $S_{vfp}^s$ & Susceptible sterile female vector & 0  \\
        $S_{vfp}^w$ & Susceptible pregnant female vector  & \\ & with Wolbachia & 0 \\
        \hline
        $I_{vf}$ & Infected female vector &  \\ & without Wolbachia & 1~500~000  \\
        $I_{vf}^w$ & Infected female vector &  \\ & with Wolbachia & 0 \\
        $I_{vfp}$ & Infected pregnant female vector  &  \\ & without Wolbachia & 500~000  \\
        $I_{vfp}^s$ & Infected sterile female vector & 0   \\
        $I_{vfp}^w$ & Infected pregnant female vector  & \\ & with Wolbachia & 0 \\
        \hline
        $A$ & Aquatic stage vector  & \\ & without Wolbachia & 25~000~000  \\
        $A^w$ & Aquatic stage vector  & \\ & with Wolbachia & 0 \\
        \hline
    \end{tabular}
    \caption{Compartments in the dengue transmission model with Wolbachia. The initial values for the vector population are divided based on the total dengue susceptible and infected vectors from \cite{escaner}.}
    \label{compartmentnames}
\end{table}

\begin{table}
    \centering
    \begin{tabular}{|r|l|c|c|}
        \hline
        Parameter & Description (units) & Value & Reference  \\
        \hline
        $b_h$ & birth rate for humans (1/day) & 0.00085/7 & \cite{escaner} \\
        $\mu_h$ & death rate for humans (1/day) & 0.00045/7 & \cite{escaner} \\
        $\alpha$ & healthcare-seeking fraction & 0.2 & \cite{escaner} \\
        $\gamma$ & nonhealthcare-seeking recovery rate (1/day) & 0.5/7 & \cite{escaner} \\
        $\theta$ & healthcare-seeking recovery rate (1/day) & 1/7 & \cite{escaner} \\
        \hline
        \hline
        $B$ & mosquito biting rate (1/day) & 1/7 & \cite{escaner} \\
        $C_{hv}$ & dengue transmission probability from & & \\ & human to mosquito & 0.75 & \cite{escaner} \\
        $C_{vh}$ & dengue transmission probability from & & \\ & mosquito (without Wolbachia) to human & 0.375 & \cite{escaner} \\
        $C_{vh}^w$ & dengue transmission probability from & & \\
        & mosquito (with Wolbachia) to human & 0 & Assumed \\
        \hline
        \hline
        $\sigma$ & mating rate (1/day) & 1 & \cite{qumodel} \\
        $\phi$ & egg-laying rate of female mosquito & & \\ & without Wolbachia (1/day) & 13 & \cite{qumodel} \\
        $\phi^w$ & egg-laying rate of female mosquito & & \\ & with Wolbachia (1/day) & 11 & \cite{qumodel} \\
        $v^w$ & maternal transmission fraction of Wolbachia & 0.95 & \cite{qumodel} \\
        $v$ & maternal transmission leakage fraction & 0.05 & \cite{qumodel} \\
        $\psi$ & per capita development rate of & & \\ & aquatic-stage mosquitoes (1/day) & 1/8.75 & \cite{qumodel} \\
        $b_m$ & male birth probability & 0.5 & \cite{qumodel} \\
        $b_f$ & female birth probability & 0.5 & \cite{qumodel} \\
        $\mu_{a}$ & per capita death rate of & & \\ & aquatic-stage mosquitoes & 0.02 & \cite{qumodel} \\ 
        $\mu_{f}$ & per capita death rate of Wolbachia-uninfected & & \\ & female mosquitoes (1/day) & 1/17.5 & \cite{qumodel} \\ 
        $\mu_{f}^w$ & per capita death rate of Wolbachia-infected & & \\ & female mosquitoes (1/day) & 1/15.8 & \cite{qumodel} \\
        $\mu_{m}$ & per capita death rate of Wolbachia-uninfected & & \\ & male mosquitoes (1/day) & 1/10.5 & \cite{qumodel} \\
        $\mu_{m}^w$ & per capita death rate of Wolbachia-infected & & \\ &  male mosquitoes (1/day) & 1/10.5 & \cite{qumodel} \\
        \hline
    \end{tabular}
    \caption{The values of the parameters used in the proposed dengue transmission model with Wolbachia.}
    \label{parameternames}
\end{table}

\begin{figure}[h!]
    \centering
    \includegraphics[width=\textwidth]{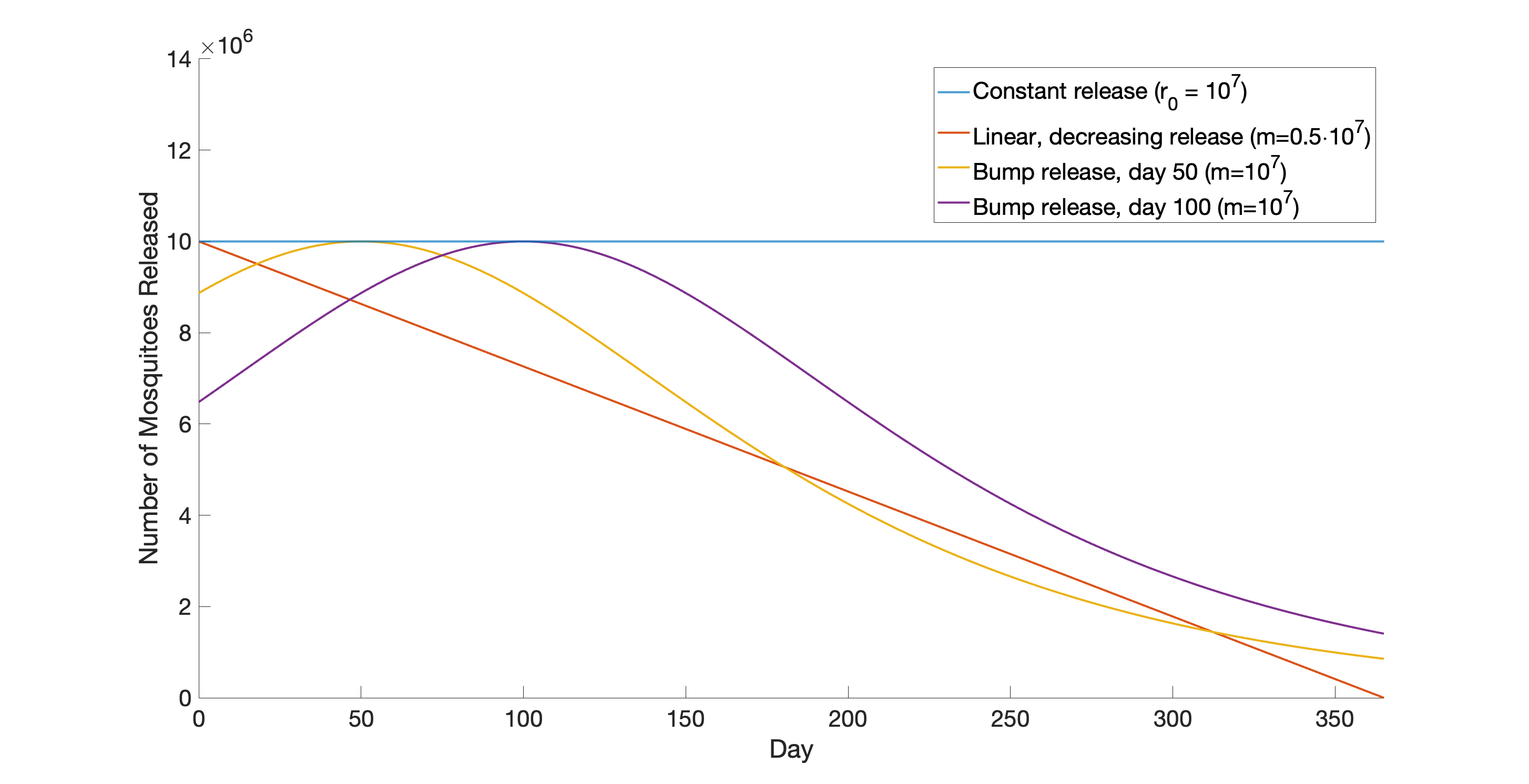}\\
    \includegraphics[width=.9\textwidth]{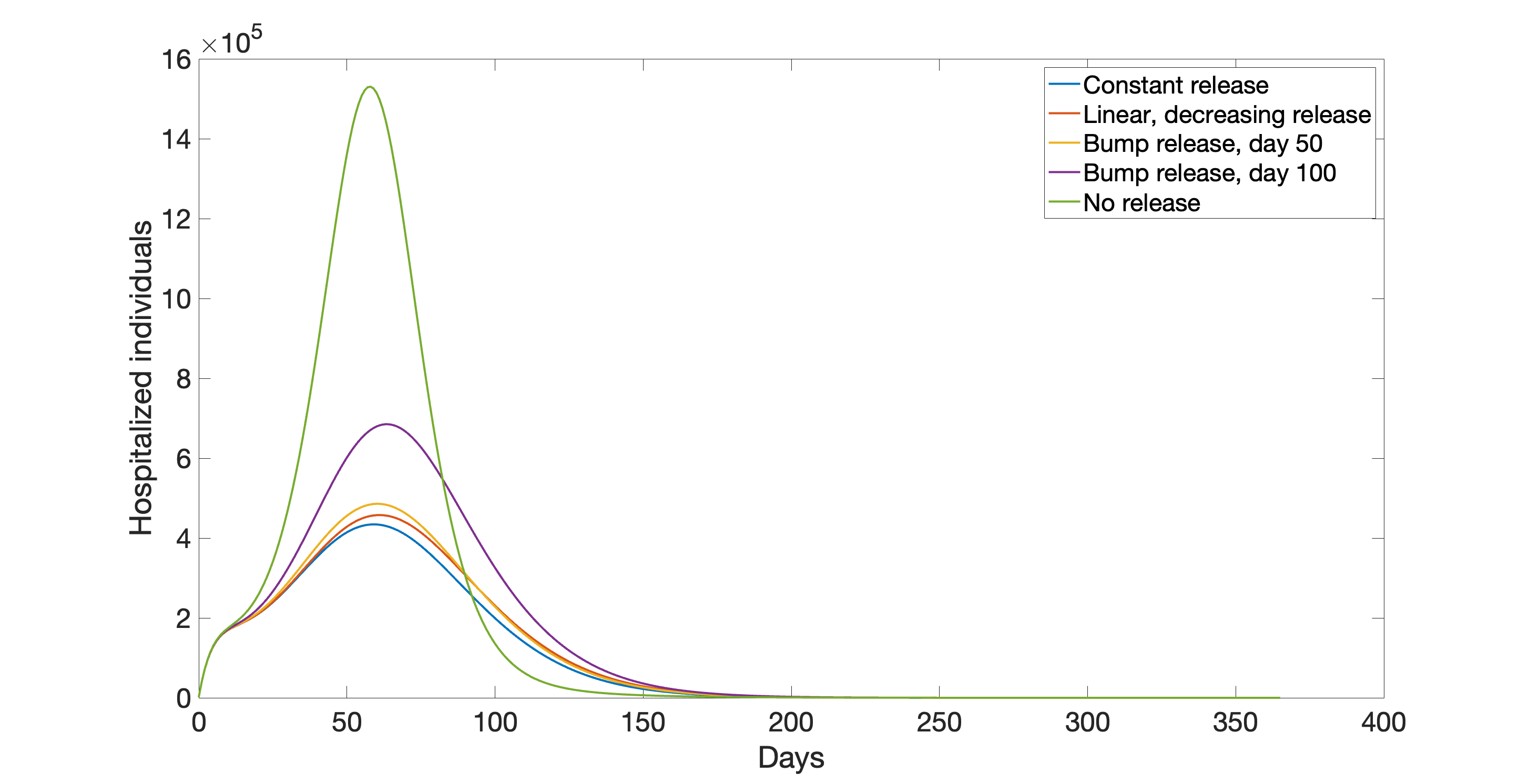} 
    \caption{Different release schemes: constant release (blue); linear and decreasing release, which simulates a peak at the start of the release period (orange); and a bump release, which simulates a peak at a later day—specifically, day 50 (yellow) or day 100 (violet). The number of hospitalized individuals under these different forms of $r(t)$ is shown for the case in which all release schemes have the same peak.}
    \label{samepeak}
\end{figure}

We investigate different forms of the release function $r(t)$ and their effects on human dengue infection prevalence.
%We explore different forms of the release function $r(t)$ and examine its effects in the prevalence of dengue infection in humans. 
The different forms of $r(t)$ are: \textit{constant}: $r(t) = r_0$; \textit{linear, decreasing}: $r(t) = -\dfrac{2m}{365} (t-365), \ m \in \mathbb{R}$; \textit{bump function with peak at day 50}: $r(t) = m \sech (0.01(t-50))$; and \textit{bump function with peak at day 100}: $r(t) = m \sech (0.01(t-100))$. A constant $r(t)$ considers a single release policy in a year. Linear and decreasing $r(t)$ simulates a release policy in which a large number of Wolbachia-infected mosquitoes is initially released into the wild, then steadily declines to zero. A bump function $r(t)$ represents a release policy in which the number of Wolbachia-infected mosquitoes released peaks at a later time and then steadily decreases. Figure \ref{samepeak} shows the graphs of the different release functions. To compare the different forms of $r(t)$, we consider both the case in which the peak of the release is the same and the case in which the total number of mosquitoes released over time is the same.

Figure \ref{samepeak} also shows the number of individuals infected with dengue, considering the same peak for $r(t)$. It shows that, given the same peak for $r(t)$, constant and linear release schemes work best, reporting a 71.6\% reduction in cases. Bump release works well in reducing dengue cases, provided that a peak release happens early. We see that a peak bump release at day 50 leads to a 30\% reduction in cases versus peak bump release at day 100. This suggests that regardless of the form of $r(t)$, an earlier peak leads to a greater reduction in dengue cases.

\begin{figure}[htb!]
    \centering
    \includegraphics[width=0.9\textwidth]{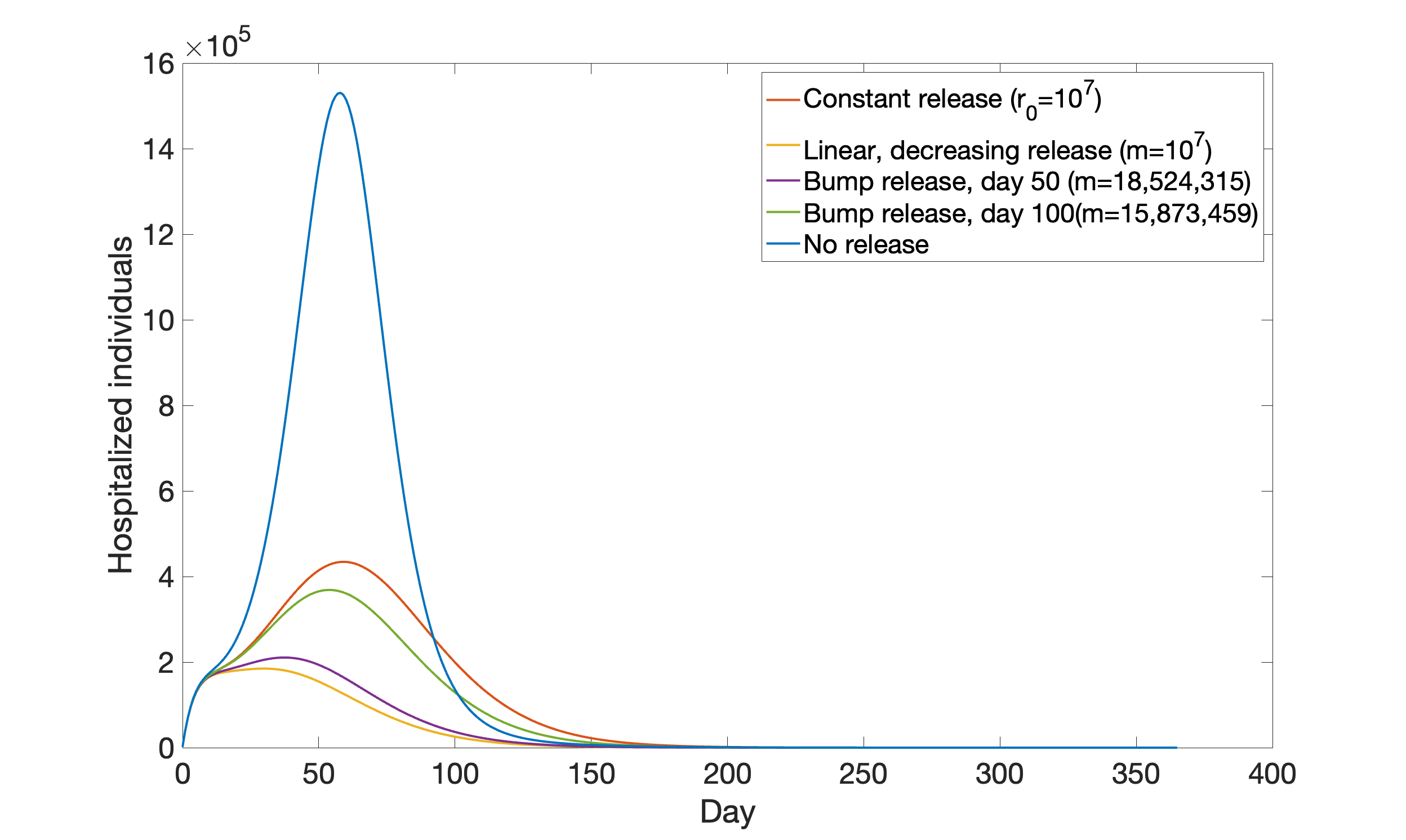}
    \caption{The number of hospitalized individuals under the different forms of $r(t)$ for the case in which all release schemes have the same total number of Wolbachia-infected larva released over time.}
    \label{sametotal}
\end{figure}

Figure \ref{sametotal} shows the number of individuals infected with dengue, considering the case in which the total number of mosquitoes released over time is the same. In this case, linear release would have the highest peak, followed by the bump releases, then the constant release. We see that despite having the same total number of Wolbachia-infected larvae released, the linear release scheme gave the highest reduction in dengue cases, with the bump release (with peak at day 50) having similar results. 

Figures \ref{samepeak} and \ref{sametotal} demonstrate that an earlier and higher peak release leads to lower prevalence of dengue in the human population. Because of this, we are interested in finding the optimal $r(t)$ that minimizes not only the number of individuals infected by dengue but also the total cost of releasing Wolbachia-infected larvae. In this work, we formulate two optimization problems, a single objective problem and a multiobjective optimization problem, which are discussed in detail in the next sections.

\subsection{Single objective optimization}\label{soo}
Single-objective optimization is first considered to obtain a clear and implementable release policy by consolidating the economic cost of Wolbachia deployment and the societal burden of dengue into a single objective function that reflects fixed budgetary priorities in public health planning.

Let $r(t)$ be the number of Wolbachia-infected mosquitoes released at time $t$ and $J(t)$ be the number of dengue healthcare-seeking individuals at time $t$. We aim to find $r^*(t)$ which minimizes the cost functional 
\[ F(r) = \sum_{t=1}^{T} [C_r r(t)+C_J J(t)].\] 
Here, $T$ is the number of days in consideration, $C_r$ is the cost of releasing Wolbachia-infected mosquitoes, and $C_J$ is the daily cost of hospitalization. 

Furthermore, we interpret $C_r r(t)$ as the cost of releasing Wolbachia-infected mosquitoes in a single day and we interpret $C_J J(t)$ as the daily economic cost of having hospitalized individuals due to dengue. Note that $C_J$ does not account only for hospitalization costs, as it takes into account the productivity losses due to hospitalization and deaths as well. 

We want $r(t)$ to be piecewise-constant, with $N$ pieces, to simulate a government policy that changes after a certain time period. With this, the release function $r(t)$ takes the form \[r(t) = \tilde{r}\left(\left\lceil \dfrac{tN}{T} \right\rceil\right),\] where $\tilde{r}(i)$ gives the release policy at the $i$th period, where $i=1,\ldots,N$.

Based on a study by Soh \cite{Cr}, the steady state cost of Wolbachia suppression programs is $40$ million SGD (Singaporean dollars). Assuming a weekly release of 7 million mosquitoes \cite{7mmosquitoes}, we get the cost (in Philippine peso) $C_r=$ PHP $4.85$ per mosquito released. Based on a study by Cheng \cite{Cj}, the total societal cost of dengue in the Philippines from 2016 to 2020 is PHP $76,275,822,554$. Dividing this by the total number of cases and considering that the computed costs is based on the fact that the average healthcare for dengue lasts for seven days, we get a daily cost of hospitalization $C_J=$ PHP $3,401.52$ per individual.

We consider two constraints for the problem. First, we consider the budget $\tilde{B}$ of the government in releasing Wolbachia-infected mosquitoes. We express this as $\displaystyle\sum_{t=1}^T C_r r(t) \leq \tilde{B}.$ Next, we consider a maximum production capacity $P$ for Wolbachia-infected mosquitoes. We want $P$ to be time-varying to simulate an improving production capacity over time and so, we have $r(t) \leq P(t),$ where $P:[1,T] \to \mathbb{R}_+$ is a bounded function. 

Thus, we have the following optimization problem:
\begin{align}
    \begin{split}
        \min_r & \quad \displaystyle\sum_{t=1}^{T} \ [C_r r(t)+C_J J(t)], \\
\text{subject to} & \ 0\leq \displaystyle\sum_{t=1}^T C_r r(t) \leq \tilde{B}, \\
& \ 0 \leq r(t) \leq P(t).
    \end{split}
    \label{proposedconstrainedproblem}
\end{align}
To establish the existence of solutions to Problem \eqref{proposedconstrainedproblem}, we rely on the form in consideration for $r(t)$ being piecewise-constant. This allows the conversion of the given problem into a constrained optimization problem in $\mathbb{R}^N$. 

Since $r(t)$ is piecewise-constant with $N$ pieces, it follows that the cost of releasing Wolbachia-infected mosquitoes over a period of $\ell = \left\lceil \dfrac{T}{N} \right\rceil$ days is $\ell C_r\tilde{r}_i$, where $\tilde{r}_i = r(\ell(i-1)+1), \ i=1,\ldots N$. With this, we rewrite Problem \eqref{proposedconstrainedproblem} as follows:
\begin{align}
    \begin{split}
        \min_{\tilde{r}} & \quad \sum_{i=1}^N \ell C_r \tilde{r}_i+\displaystyle\sum_{t=1}^{T} \ C_j J(t), \\
\text{subject to} & \ 0\leq \displaystyle\sum_{t=1}^T \ell C_r  \tilde{r}_i \leq \tilde{B}, \\
& \ 0 \leq \tilde{r}_i \leq P(\ell(i-1)+1).
    \end{split}
    \label{newproposedproblem}
\end{align}

\begin{theorem}
    The optimization problem \eqref{newproposedproblem} has a minimizer in \[X = \left\{ (\tilde{r}_1,\ldots,\tilde{r}_N) \in  \mathbb{R}^N_{{+}} \ \mid \ 0\leq \tilde{r}_i\leq P(\ell(i-1)+1), i=1,\ldots,N, 0\leq \sum_{i=1}^N \ell C_r \tilde{r}_i\leq \tilde{B} \right\}.\]
    \begin{proof}
        We see that $X$ is nonempty since the zero vector in $\mathbb{R}^N$ is in $X$.
        
        {Now, note that the first term of the objective function is a polynomial in $r_i$. So, it is continuous in $X$. To establish the continuity of the second term, let $F(t,x(t))=[F_1(t,x(t)),\ldots,F_{18}(t,x(t))]^T$ be the vector-valued function corresponding to the right-hand side of the equations in \eqref{modelhuman}-\eqref{modelvec2}. See that the derivative of $F$ with respect to the function $r:=r(t)$ is given by \[\dfrac{dF}{dr} = [0,\ldots,0,1]^T,\] which is continuous. So, $F$ is continuous with respect to the function $r$, and consequently, $J(t)$ is continuously dependent on $r$. With this, it follows that the second term is continuous on $X$.} So, the objective function is continuous in $X$.
        
        To prove compactness of $X$, we argue that $X$ is closed and bounded. By definition, $X\subseteq [0,M]^N$, with $M:=\max P(t)$, which is a bounded set. It follows that from this that $X$ is bounded. Now, we can rewrite $X$ as follows: \[X = \tilde{X} \  \cap \ \left\{(r_1,\ldots,r_N)\in \mathbb{R}^N_{{+}} \ | \ 0\leq \sum_{i=1}^N r_i \leq \tilde{B}\left(\ell C_r \right)^{-1} \right\}, \] where \[\tilde{X}=[0,P(1)] \times [0,P(\ell+1)] \times \cdots \times [0,P(\ell(N-1)+1)]. \] Since each interval of the form $[0,P(\ell(i-1)+1)]$ is closed, it follows that $\tilde{X}$ is closed. This means that $X$ is the intersection of a closed set and a closed upper half-space {bounded} by the hyperplane $\displaystyle\sum_{i=1}^N r_i = \tilde{B}\left(\ell C_r\right)^{-1}$. Hence, $X$ is closed.
        
        Since $X$ is closed and bounded, it follows that $X$ is compact, and thus, by Weierstrass theorem, the objective function admits a minimizer in $X$.
      
    \end{proof}
\end{theorem}

To find the optimal solution $r^*$, we start by solving the system of equations in \eqref{modelhuman}-\eqref{modelvec2} as in Section \ref{simulation}. Then, we proceed by employing an interior-point method to improve $r(t)$. We do this process iteratively until we obtain $r^*$. We utilize the Matlab built-in function \texttt{fmincon}, which makes use of a trust-region method based on the interior-point method to implement this optimization technique \cite{fmincon}.

With this single objective optimization problem, equal priority is given to the budget constraint of releasing Wolbachia and the cost of hospitalization. In the next section, we introduce a multiobjective optimization to explicitly characterize the trade-off between intervention expenditure and dengue burden, allowing the identification of Pareto-optimal release strategies without imposing a priori weights on competing costs.

\subsection{Multiobjective optimization}

We consider a multiobjective optimization problem, wherein we minimize the total cost of releasing Wolbachia-infected mosquitoes $\sum_{t=1}^T C_r r(t)$ and the total societal cost $\sum_{t=1}^T C_J J(t)$, simultaneously. The multiobjective optimization problem is written as follows:

\begin{align}
&\min_{r(t)} 
\begin{cases}
	\sum\limits_{t=1}^T \ C_J J(t)
    \label{eq:optProb1}\\
	\sum\limits_{t=1}^T  \ C_r r(t)  \\
\end{cases}\\
&\text{subject to } r(t) \leq P(t).\label{eq:constraint2}
\end{align}

To solve this problem, we employ the $\varepsilon-$constraint method, where we use one objective function as a constraint and change its upper bound to produce a set of optimal solutions for the problem, called the Pareto-optimal front \cite{pareto}. In this case, we let the societal cost be the objective function and the total cost of releasing Wolbachia-infected mosquitoes be a constraint. We consider $K=100$ upper bounds for this, ranging from 0 to PHP $500,000,000$. The optimization problem \eqref{eq:optProb1}-\eqref{eq:constraint2} can be rewritten as follows:

\begin{subequations}
\begin{alignat*}{2}
&\!\min_{r(t)}        &\qquad \sum_{t=1}^T C_j J(t) \\
&\text{subject to} &       \sum_{t=1}^T C_r r(t) & \leq \tilde{B}_k, \quad {k = 1,\ldots,K},\\
&                  &       r(t) &\leq P(t).\label{eq:constraint2}
\end{alignat*}
\end{subequations}

Since the optimization is converted into a single objective optimization problem, we employ the same method described in Section \ref{soo} to numerically solve the problem.

\section{Results and discussion}
\label{sec:results}
In this section, we present the results of the two optimization problems described in the previous section. We begin with the single-objective problem, examining the effects of varying the production capacity and the cost of releasing Wolbachia-infected mosquitoes, as shown in subsections \ref{subsec:maxprod} and \ref{subsec:reducedprice}, respectively. Subsection \ref{sec:multiresult} presents the results of the multiobjective problem.

In our simulations, we use the initial values from Escaner and de  los Reyes V\cite{escaner}, but scaled to a smaller population (i.e., from an initial susceptible population of 50 million to 2960000, the population of Quezon City). This means that if $X$ is the vector of initial values, then we consider a new initial value of 
\[ X_{new} = \dfrac{2960000}{50000000}X.\] 
We consider a time period of $T=$365 days for the simulations.
In addition, we consider $r(t)$ as a piecewise constant function. We assume that $r(t)$ consists of $N=12$ pieces to simulate a monthly change in policy.

\begin{figure}[htb!]
	\centering
	\includegraphics[width=.9\textwidth]{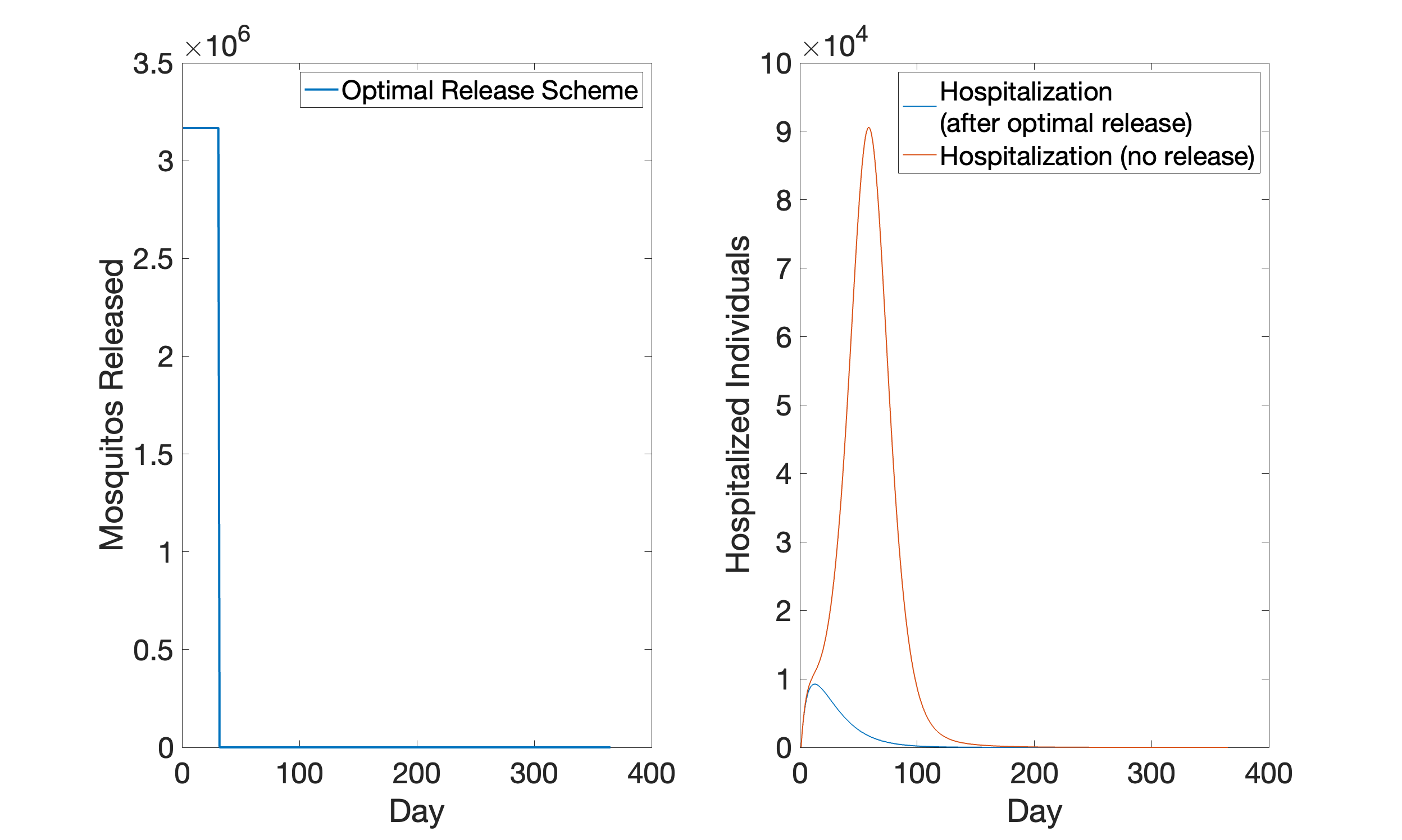}
	\caption{Optimal release schedule and its effect on healthcare-seeking individuals}
	\label{noupperboundrelease}
\end{figure}

In Figure \ref{noupperboundrelease}, we see the optimal release of Wolbachia under the constraint that $r(t)\geq 0$. This result illustrates that a higher initial release is key to minimizing hospitalizations. The obtained $r(t)$ suggests that we release 3.17 million Wolbachia-infected mosquitoes initially then drop to 18 mosquitoes daily by day 32. Doing this would lead to a peak hospitalization of 9561 people at day 14, which is a 90\% reduction compared to not releasing Wolbachia-infected mosquitoes at all. Under this release schedule, the total cost is PHP $1,668,558,751.79$, where PHP $476,155,281.37$ is the cost from the release of Wolbachia-infected mosquitoes and PHP $1,192,403,470.42$ is the societal cost.

\subsection{Effects of Maximum Production Capacity}
\label{subsec:maxprod}
While the previous result is ideal in the sense that it yields the lowest objective value, we must consider that maximum production is not typically achieved during the initial days of operation. Therefore, we introduce a maximum production $P(t)$ for the release of Wolbachia-infected mosquitoes. Specifically, we assume that $P(t)$ increases linearly until it reaches a peak value of 3.5 million on Day 94 (Month 3).

\begin{figure}[htb!]
	\centering
	\includegraphics[width = 0.8\textwidth]{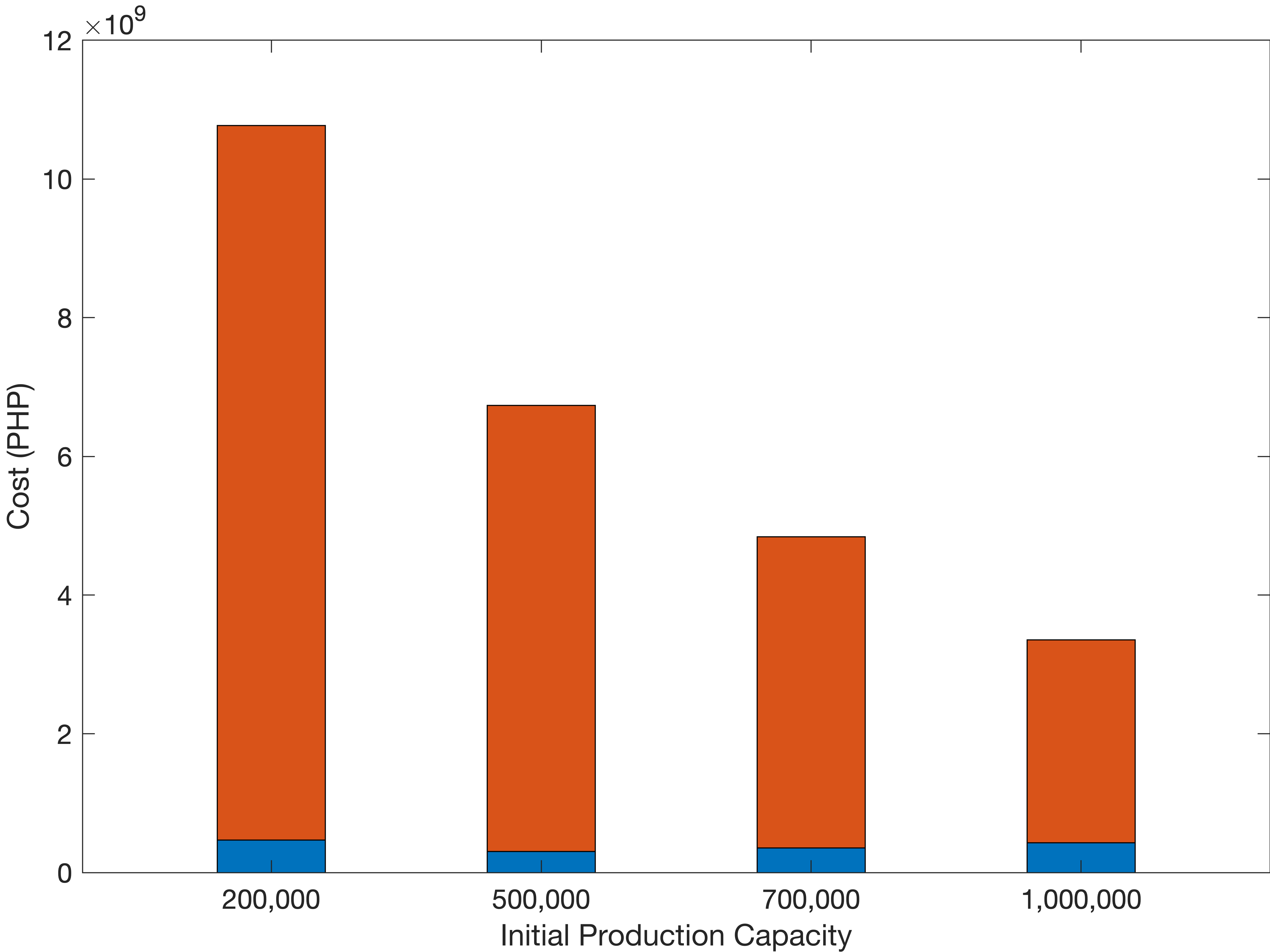}
	      \caption{Total annual cost, including the annual cost of releasing Wolbachia-infected mosquitoes (blue) and the annual societal cost due to dengue infections (red), under different initial production capacities.}
	\label{compinitcapacity}
\end{figure}

In Figure \ref{compinitcapacity}, we see that a higher initial capacity reduces the societal cost, and so, it also reduces the total cost. However, since the societal cost is not necessarily shouldered by the government, we focus on the cost of releasing Wolbachia-infected mosquitoes. For initial maximum production capacity of 500,000 to 1,000,000, we see that a higher cost of releasing Wolbachia-infected mosquitoes leads to a lower societal cost. For the case where the initial maximum production capacity is 200,000, we see both a higher societal cost and a higher cost of releasing Wolbachia-infected mosquitoes. This suggests that a high initial maximum production capacity is important for an effective release scheme of Wolbachia-infected mosquitoes.

\begin{figure}[htb!]
	\centering
	\includegraphics[width = \textwidth]{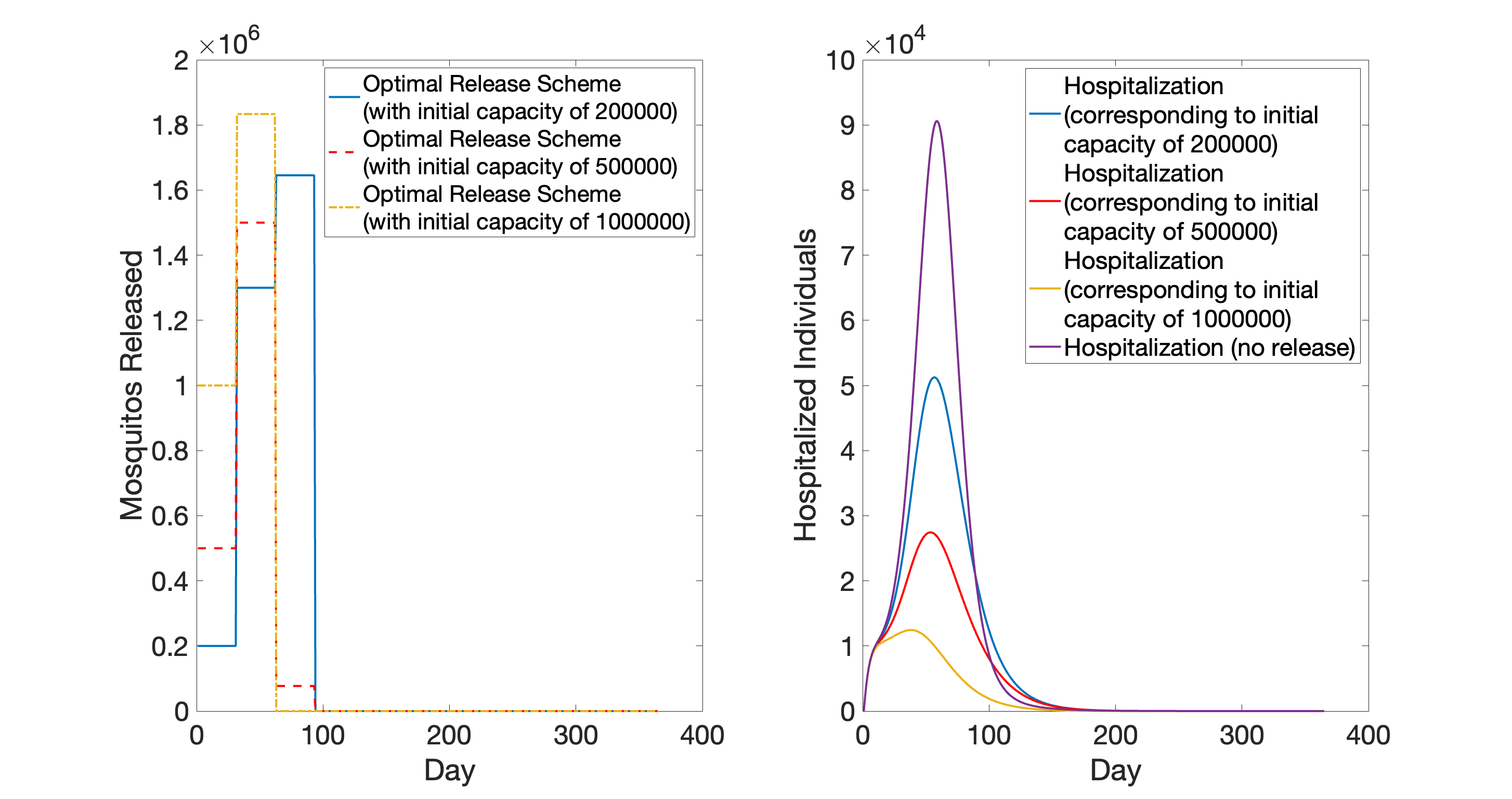}
	\caption{Optimal release schedule considering different initial capacities. The left plot shows the optimal release schedule under different initial capacities while the right plot shows the number of hospitalized individuals at time $t$ under the different release schemes.}
	\label{diffinitcapacity}
\end{figure}

Figure \ref{diffinitcapacity} shows the optimal release schedule given different initial production capacities. If we look at the hospitalizations due to dengue, we see that a higher initial capacity leads to a lower hospitalization, and hence, lower societal cost. If we consider an initial capacity of 700,000, we see a 79\% decrease in hospitalization, and if we consider an initial capacity of 1,000,000, we see an 86\% decrease in hospitalization.

Now, we observe the effect in hospitalizations of an earlier peak or a later peak in production of Wolbachia-infected mosquitoes. For our simulation, we use an initial capacities of 1 million and 500,000.

\begin{figure}[htb!]
	\centering
	\includegraphics[width = .8\textwidth]{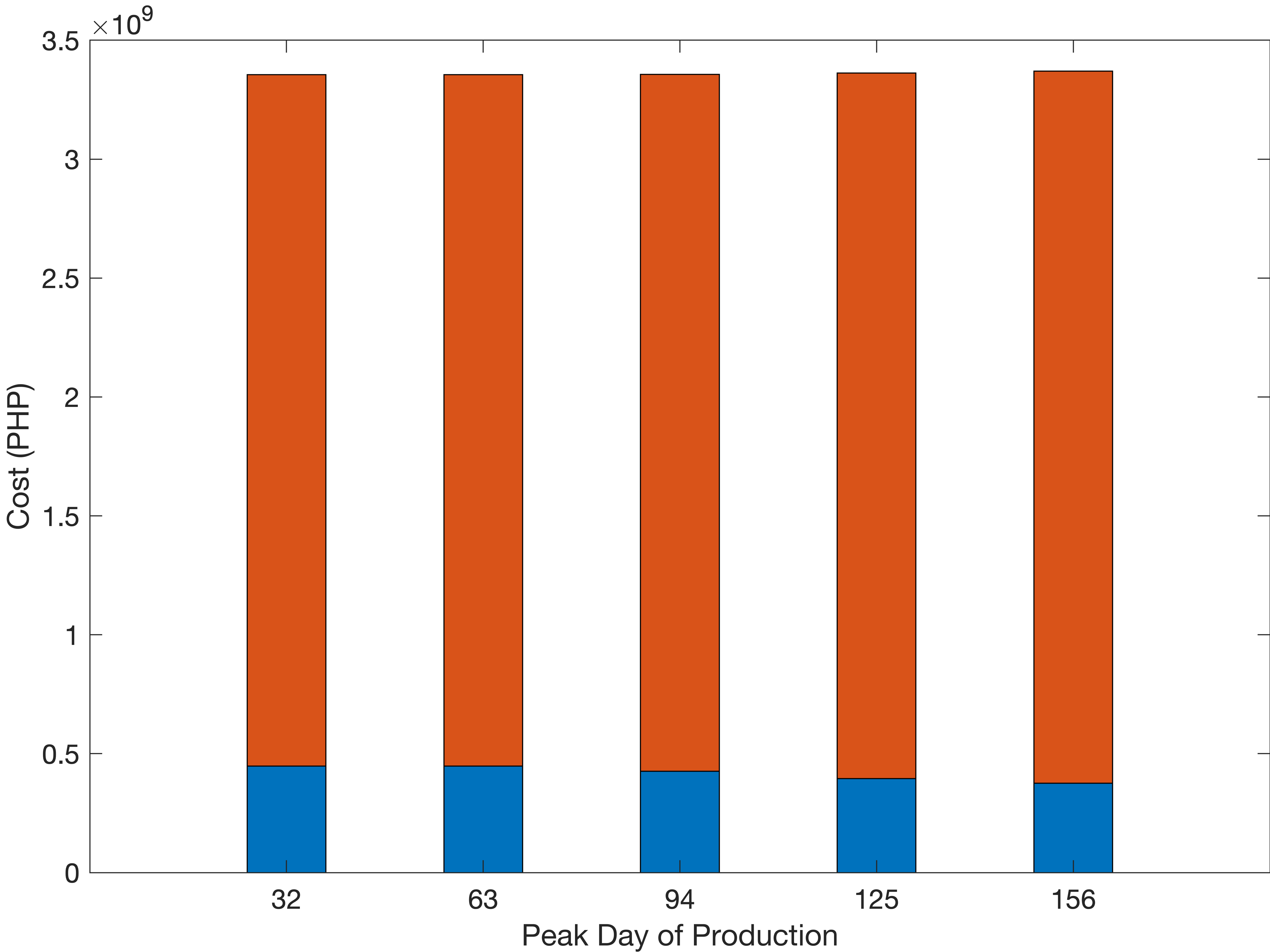}
        \caption{Total annual cost, including the annual cost of releasing Wolbachia-infected mosquitoes (blue) and the annual societal cost due to dengue infections (red), under different peak days of production, given initial capacity of 1 million.}
	\label{comppeakday1m}
\end{figure}

Figure \ref{comppeakday1m} shows the cost of releasing Wolbachia-infected mosquitoes and the societal cost under different peak days of production capacity, given an initial capacity of 1 million. 
We observe that if the production facility reaches maximum output early, the societal cost due to dengue decreases slightly. However, this also results in a higher cost for releasing Wolbachia-infected mosquitoes, since the optimal release tends to coincide with the maximum possible production. Overall, the total cost does not change significantly, regardless of when the production facility reaches its maximum output.

\begin{figure}[htb!]
    \centering
    \includegraphics[width=\textwidth]{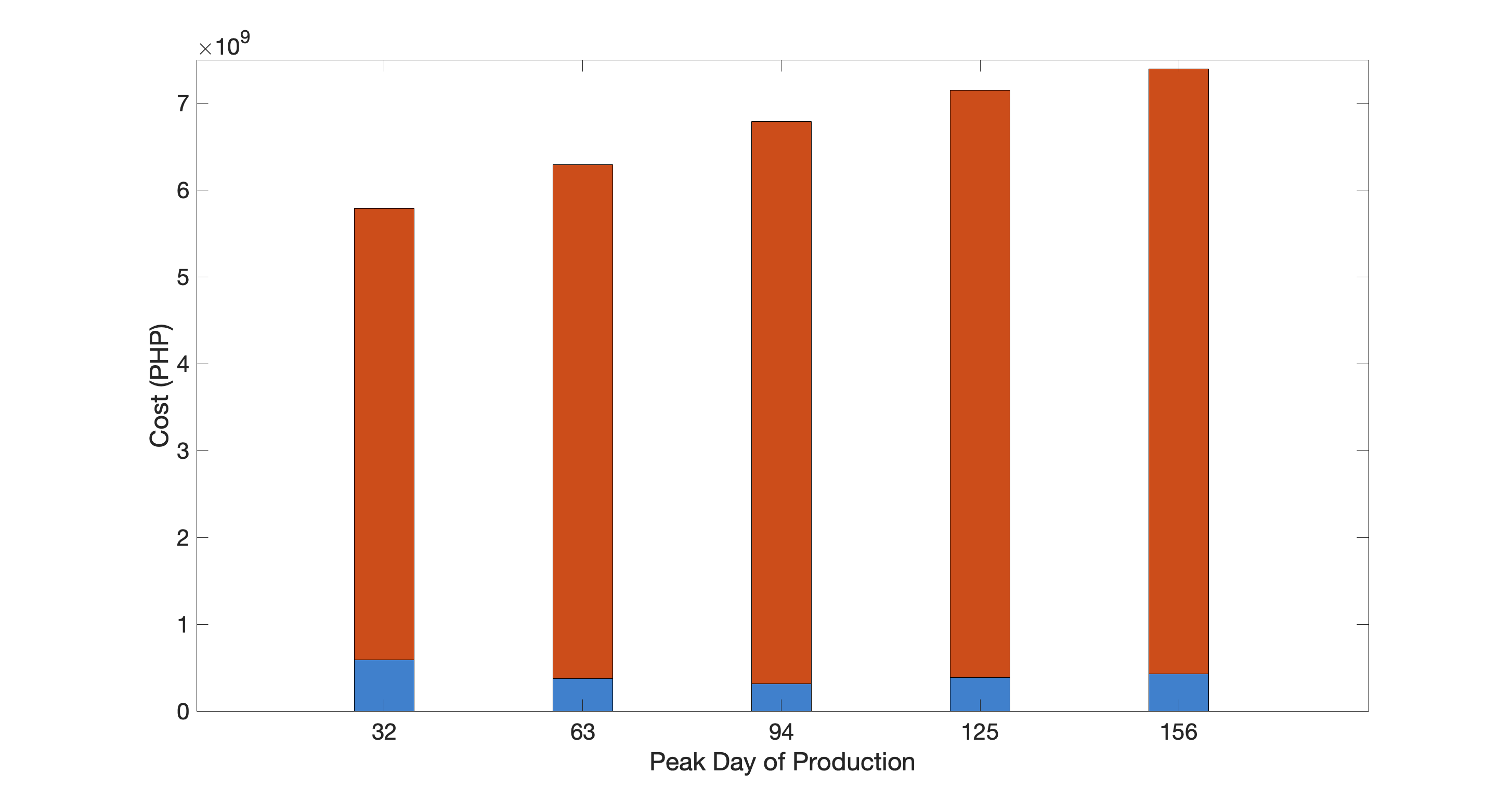}
    \caption{Total annual cost, including the annual cost of releasing Wolbachia-infected mosquitoes (blue) and the annual societal cost due to dengue infections (red), under different peak days of production, given initial capacity of 500,000.}
    \label{comppeakday500k}
\end{figure}

\begin{figure}[htb!]
    \centering
    \includegraphics[width=\linewidth]{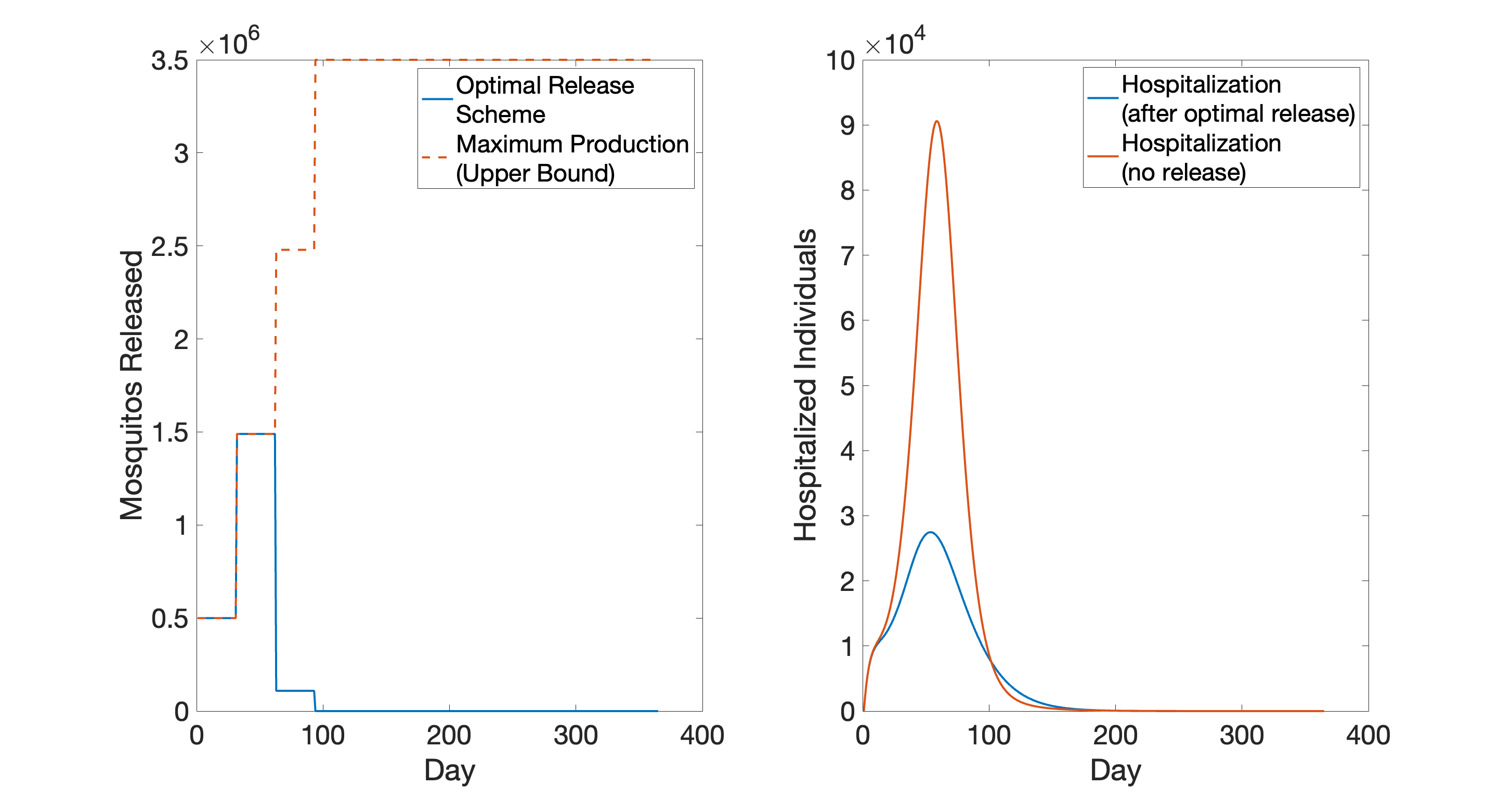}
    \caption{Optimal release schedule considering an initial capacity of 500,000 and a peak of production starting on Day 94.}
    \label{peakday94}
\end{figure}

\begin{figure}[htb!]
    \centering
    \includegraphics[width=\linewidth]{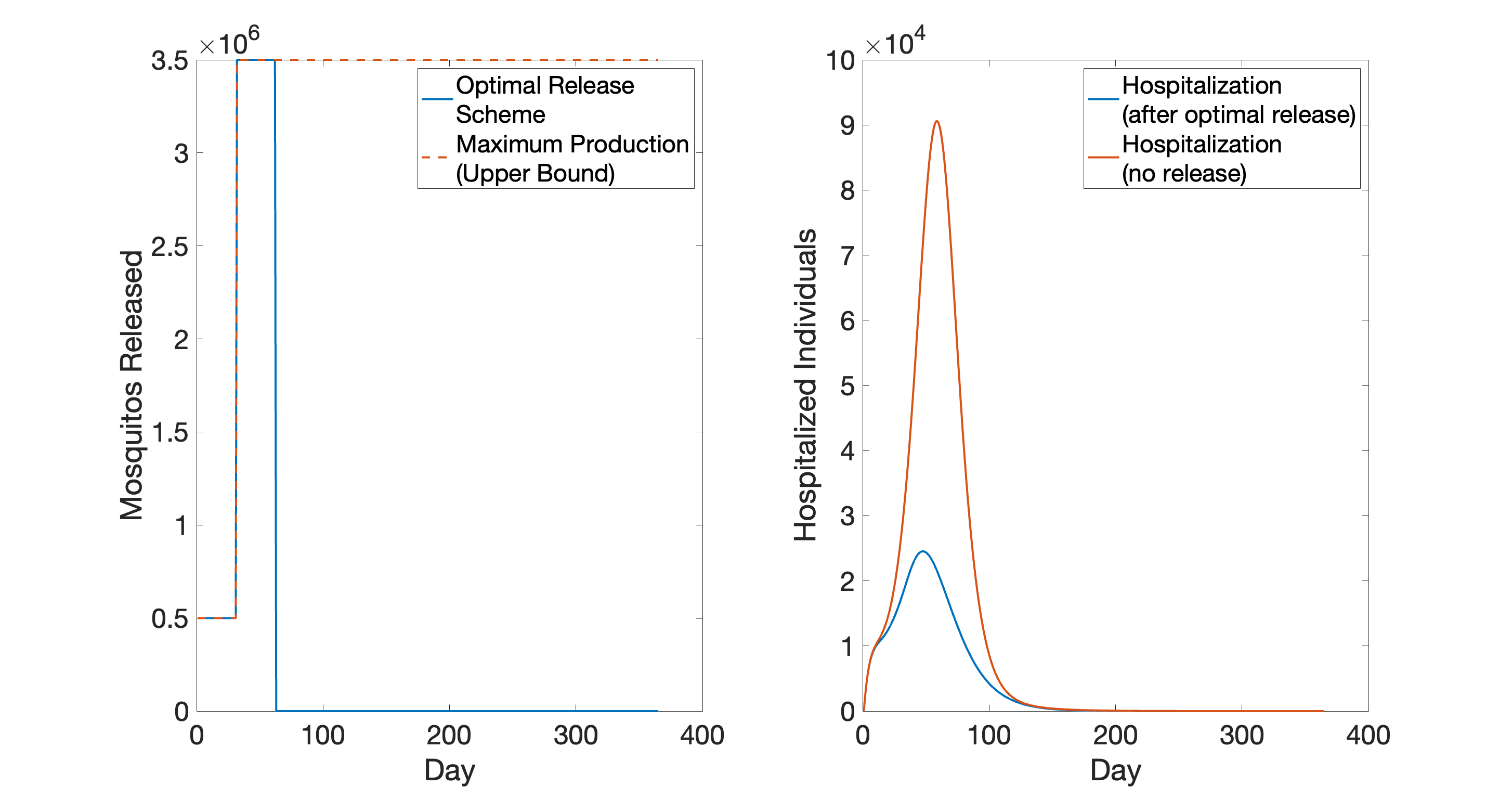}
    \caption{Optimal release schedule considering an initial capacity of 500,000 and a peak of production starting on Day 32.}
    \label{peakday32}
\end{figure}

Figure \ref{comppeakday500k} shows the cost of releasing Wolbachia-infected mosquitoes and the societal cost under different peak days of production capacity, given an initial capacity of 500,000. Unlike in the case where the initial capacity is 1 million, an increase in the total cost is seen as we delay the peak day of production. In Figures \ref{peakday94} and \ref{peakday32}, we see that the peak hospitalization is lower if production capacity peaks early. This suggests that for lower initial capacities, an earlier peak day of production would lead to a lower societal cost, and hence, a lower total cost associated with dengue infections.

\subsection{Effects of a Lower Price of Releasing Wolbachia-infected mosquitoes}
\label{subsec:reducedprice}
However, note that the cost of labor in the Philippines is lower than in Singapore. So, we want to determine the effects of a lower unit price for Wolbachia-infected mosquitoes to the number of hospitalizations due to dengue. For our simulations, we assume an initial maximum production of 1,000,000 mosquitoes and a peak maximum production of 3,500,000 mosquitoes achieved at the third month. Table \ref{tab:diffunitprices} shows the total release of Wolbachia-infected mosquitoes and {total societal cost} given different unit prices.

\begin{table}[htb!]
	\centering
	\begin{tabular}{|c|c|c|}
		\hline
		Unit price (Php) & Total release of mosquitoes & Total societal cost (Php) \\
		\hline
		4.85 & 87,833,750 & 2,930,247,072.05 \\
		\hline
		4.00 & 87,833,559 & 2,930,247,159.25 \\
		\hline
		3.00 & 87,833,390 & 2,930,247,217.25 \\
		\hline
		2.00 & 87,833,333 & 2,930,247,240.27 \\
		\hline
	\end{tabular}
	\caption{Comparison of total release of mosquitoes and total societal cost under different unit prices for Wolbachia-infected mosquitoes (given an initial capacity of 1 million).}
    \label{tab:diffunitprices}
\end{table}

We see that as the unit price decreases, the total release of mosquitoes decreases, although the decrease is not significant. Furthermore, reducing the unit price from PHP $4.85$ to PHP $2.00$ (58.8\% decrease in unit price) leads to a 0.0005\% decrease in the total release of mosquitoes. This suggests that, given a high initial capacity, the unit price of Wolbachia-infected mosquitoes does not significantly affect the amount of mosquitoes to be released. 

\begin{table}[H]
	\centering
	\begin{tabular}{|c|c|c|}
		\hline
		Unit price (Php) & Total release of mosquitoes & Total societal cost (Php) \\
		\hline
		4.85 & 64,530,150 & 6,465,531,657.34 \\
		\hline
		4.00 & 74,393,853 & 6,422,124,708.67 \\
		\hline
		3.00 & 89,767,915 & 6,368,727,268.25 \\
		\hline
		2.00 & 111,906,046 & 6,314,028,548.29 \\
		\hline
	\end{tabular}
	\caption{Comparison of total release of mosquitoes under different unit prices for Wolbachia-infected mosquitoes (given an initial capacity of 500,000).}
    \label{tab:diffunitprices500k}
\end{table}

{Table \ref{tab:diffunitprices500k} shows the total release of Wolbachia-infected mosquitoes and the corresponding total societal cost given different unit prices, assuming an initial capacity of 500,000. In contrast to the case where the initial capacity is 1 million, a decrease in unit price leads to higher total release of mosquitoes. This suggests that a lower unit price allows for a higher level of release of Wolbachia-infected mosquitoes, leading to lower total societal costs.}

\begin{table}[H]
	\centering
	\begin{tabular}{|c|c|c|}
		\hline
		Unit price (Php) & Initial capacity of 1 million & Initial capacity of 500,000 \\
		\hline
		4.85 & 3,356,240,760.44 & 6,778,502,888.77 \\
		\hline
		4.00 & 3,281,581,394.10 & 6,719,700,122.69 \\
		\hline
		3.00 & 3,193,747,388.65 & 6,638,031,013.43 \\
		\hline
		2.00 & 3,105,913,907.17 & 6,537,840,641.12 \\
		\hline
	\end{tabular}
	\caption{{Comparison of the total cost of releasing Wolbachia-infected mosquitoes and of hospitalizations under different initial capacities and different unit prices.}}
    \label{tab:totalcost}
\end{table}

{Table \ref{tab:totalcost} shows the total cost of releasing Wolbachia-infected mosquitoes and of hospitalizations under different initial capacities and different unit prices. In both cases, the total cost decreases as unit price increases. This is evident with the decreasing total societal cost as unit price decreases, as well as the lower cost in releasing Wolbachia-infected mosquitoes as a consequence of decreasing the unit price.}

\subsection{Multiobjective optimization}\label{sec:multiresult}

\begin{figure}
    \centering
    \includegraphics[width=\linewidth]{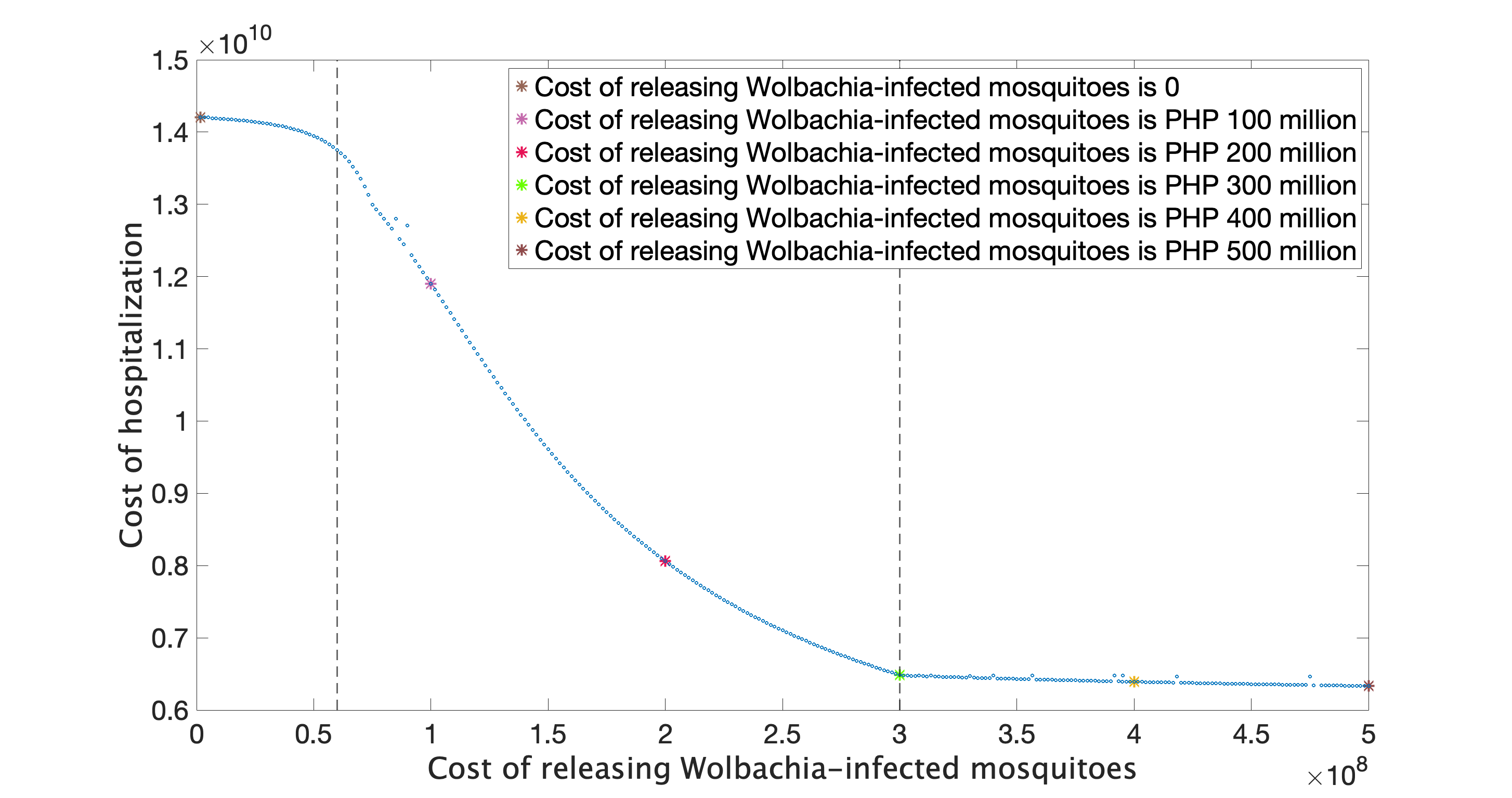} \\
    \includegraphics[width=\linewidth]{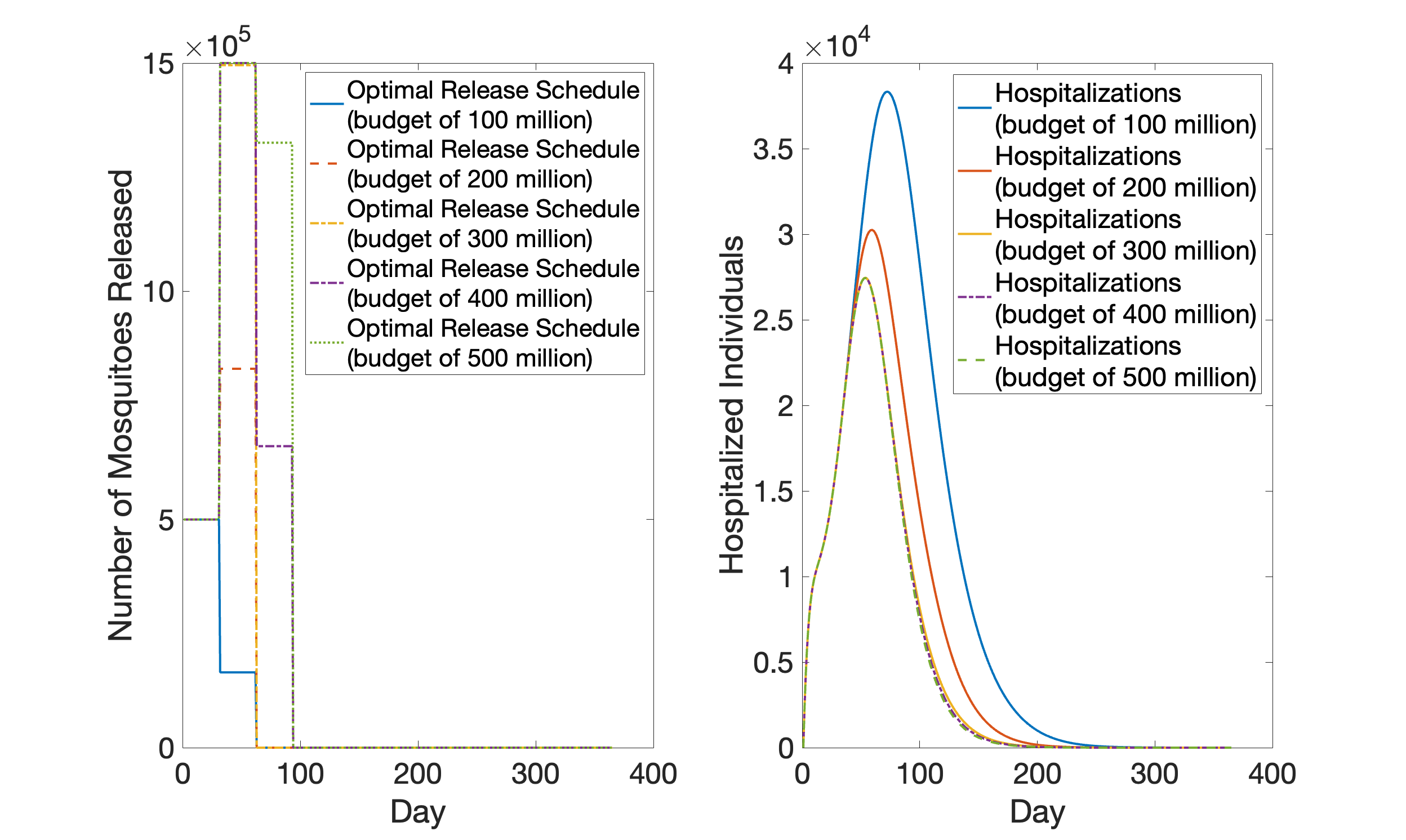}
    \caption{Pareto-optimal front for the multiobjective optimization problem. The optimal release schedule and number of hospitalized individuals for the highlighted optimal solutions are shown below the Pareto-optimal front.}
    \label{pareto}
\end{figure}

Figure \ref{pareto} shows the set of optimal solutions to the multiobjective optimization problem. We see that the societal cost remains at around PHP $14$ billion for low total release costs. Once the release cost reaches PHP $75$ million, the societal cost rapidly drops, then stabilizes at around PHP $6.5$ billion pesos once a release cost of PHP $300$ million is reached. This shows that around a release cost of PHP $75$ million to PHP $300$ million, there is a good trade-off in increasing the release of Wolbachia-infected mosquitoes. Beyond these values, a higher release leads to diminishing returns in societal cost. The multiobjective optimization approach provides a range of optimal solutions, allowing policymakers to choose an implementation strategy based on budget constraints and desired societal outcomes.

In summary, numerical simulations show that an optimal release schedule, assuming the absence of upper bound, peaks at the first period of release then drops to zero afterwards. If an upper bound is applied, then the optimal release schedule follows the trend of the upper bound initially then drops to zero after some time. This suggests the importance of a high early release of Wolbachia-infected mosquitoes. The observed trend is not affected by the cost of releasing Wolbachia-infected mosquitoes, which suggests that lowering the unit cost of Wolbachia-infected mosquitoes does not entail an increase in the amount of release. Furthermore, it is shown that the cost of the optimal release and the cost of hospitalization associated with the optimal release has an inverse relationship. With this, a multiobjective approach is applied, giving a set of optimal solutions to the problem.

\section{Conclusion}
\label{sec:conc}

In this work, we proposed a model that incorporates the Wolbachia population dynamics into the dengue transmission model. We established the existence, uniqueness, and positivity of solutions for the proposed model. Furthermore, preliminary simulations considering different release schemes motivated the formulation of the optimization problems. 
A single-objective optimization problem is considered, in which the combination of the total cost of releasing Wolbachia-infected mosquitoes and the annual societal cost is minimized.
Simulations indicate that the optimal release scheme tends to utilize the maximum production capacity at the beginning. Under the best scenario, releasing Wolbachia-infected mosquitoes could lead to a 90\% reduction in the number of hospitalized individuals. The simulations also show that the total cost of releasing Wolbachia-infected mosquitoes and the annual societal cost are inversely related.

A multiobjective optimization problem is considered, in which both the total cost of releasing Wolbachia-infected mosquitoes and the annual societal cost are minimized.
Simulations show that increasing the budget for releasing Wolbachia-infected mosquitoes initially causes a significant decrease in the societal cost, but once a certain threshold for the budget is reached, the decrease in societal cost slows down. 
The results from the multiobjective optimization problem could potentially inform government considerations on the use of Wolbachia-infected mosquitoes as a dengue suppression strategy. This approach might allow for more flexible, targeted, and cost-effective decision-making in dengue control.

It is important to note that all results presented in this study are based on theoretical analysis and numerical simulations. At present, Wolbachia-based interventions have not yet been implemented in the Philippines, and thus no empirical local data are available for validation. Nevertheless, the proposed framework provides a quantitative basis that can assist policymakers in the Philippines in planning and evaluating potential strategies should the country decide to adopt Wolbachia interventions in the future. Moreover, while the model parameters and analysis are motivated by the Philippine context, the structure of the model and the optimization framework are sufficiently general and can be adapted to other dengue-endemic countries with appropriate parameter calibration.

For future work, the effects of dengue seasonality on the release of Wolbachia-infected mosquitoes can be considered. The model can be extended to consider other suppression strategies such as pesticides and fogging. Since Wolbachia increases the prevalence of lymphatic filariasis \cite{filariasisx}, a model that looks into the transmission of both dengue and lymphatic filariasis with the presence of Wolbachia is another work that can also be pursued. Lastly, a cost-benefit analysis on the release of Wolbachia-infected mosquitoes can be done. The results of the analysis can be applied to formulate a more accurate policy proposal for the intervention.

\section*{Acknowledgment}
The authors acknowledge the \textit{DOST's Grant to Outstanding Achievements in Science and Technology} for funding support.
\newpage
\appendix \section{Appendix. Jacobian Matrix of $F$ in Theorem \ref{existence}}
\label{JacobianF}

The Jacobian matrix of $F$ in Theorem \ref{existence} is given by the block matrix 
\begin{equation}
    J = \begin{bmatrix}
        J_1 & J_2 & J_3 & J_4 & J_5 & J_6
    \end{bmatrix},
\end{equation}
where
\allowdisplaybreaks{\begin{equation*}
    \begin{split}
     J_1=   \tiny
    \begin{bmatrix}
        b_h - \mu_h + \dfrac{(BC_{vh}(I_{vfp}+I_{vf})+BC_{vwh}(I_{vfp}^s+I_{vfp}^w+I_{vf}^w))S_h}{N_h^2}-\dfrac{BC_{vh}(I_{vfp}+I_{vf})+BC_{vwh}(I_{vfp}^s+I_{vfp}^w+I_{vf}^w)}{N_h} \\
        \dfrac{(1-\alpha)(BC_{vh}(I_{vfp}+I_{vf})+BC_{vwh}(I_{vfp}^s+I_{vfp}^w+I_{vf}^w))(I_h+J_h+R_h)}{N_h^2} \\
        \dfrac{\alpha(BC_{vh}(I_{vfp}+I_{vf})+BC_{vwh}(I_{vfp}^s+I_{vfp}^w+I_{vf}^w))(I_h+J_h+R_h)}{N_h^2} \\ 0 \\
        \dfrac{BC_{hv}I_hS_{vf}^w}{N_h^2} \\ - \dfrac{BC_{hv}I_hS_{vf}^w}{N_h^2} \\
        \dfrac{BC_{hv}I_hS_{vf}}{N_h^2} \\ -\dfrac{BC_{hv}I_hS_{vf}}{N_h^2} \\ 0 \\ 0 \\
        \dfrac{BC_{hv}I_hS_{vfp}^w}{N_h^2} \\ -\dfrac{BC_{hv}I_hS_{vfp}^w}{N_h^2} \\ \dfrac{BC_{hv}I_hS_{vfp}^s}{N_h^2} \\ -\dfrac{BC_{hv}I_hS_{vfp}^s}{N_h^2} \\ \dfrac{BC_{hv}I_hS_{vfp}}{N_h^2} \\ -\dfrac{BC_{hv}I_hS_{vfp}}{N_h^2} \\ 0 \\ 0
    \end{bmatrix},
    \end{split}
\end{equation*}}
\allowdisplaybreaks{\begin{equation*}
    \begin{split}
        J_2 = \tiny
        \begin{bmatrix}
            b_h+ \dfrac{(BC_{vh}(I_{vfp}+I_{vf})+BC_{vwh}(I_{vfp}^s+I_{vfp}^w+I_{vf}^w))S_h}{N_h^2} \\
            -\gamma - \mu_h - \dfrac{(1-\alpha)(BC_{vh}(I_{vfp}+I_{vf})+BC_{vwh}(I_{vfp}^s+I_{vfp}^w+I_{vf}^w))S_h}{N_h^2} \\ 
            \dfrac{\alpha(BC_{vh}(I_{vfp}+I_{vf})+BC_{vwh}(I_{vfp}^s+I_{vfp}^w+I_{vf}^w))S_h}{N_h^2} \\ \gamma \\
            -\dfrac{BC_{hv}(S_h+J_h+R_h)S_{vf}^w}{N_h^2} \\ \dfrac{BC_{hv}(S_h+J_h+R_h)S_{vf}^w}{N_h^2} \\
            - \dfrac{BC_{hv}(S_h+J_h+R_h)S_{vf}}{N_h^2} \\ \dfrac{BC_{hv}(S_h+J_h+R_h)S_{vf}}{N_h^2} \\ 0 \\ 0 \\
            -\dfrac{BC_{hv}(S_h+J_h+R_h)S_{vfp}^w}{N_h^2} \\ \dfrac{BC_{hv}(S_h+J_h+R_h)S_{vfp}^w}{N_h^2} \\
            -\dfrac{BC_{hv}(S_h+J_h+R_h)S_{vfp}^s}{N_h^2} \\ \dfrac{BC_{hv}(S_h+J_h+R_h)S_{vfp}^s}{N_h^2} \\
            -\dfrac{BC_{hv}(S_h+J_h+R_h)S_{vfp}}{N_h^2} \\ \dfrac{BC_{hv}(S_h+J_h+R_h)S_{vfp}}{N_h^2} \\ 0 \\ 0
        \end{bmatrix},
    \end{split}
\end{equation*}}
\allowdisplaybreaks{\begin{equation*}
    \begin{split}
        J_3 = \tiny
        \begin{bmatrix}
            b_h+ \dfrac{(BC_{vh}(I_{vfp}+I_{vf})+BC_{vwh}(I_{vfp}^s+I_{vfp}^w+I_{vf}^w))S_h}{N_h^2} \\
            - \dfrac{(1-\alpha)(BC_{vh}(I_{vfp}+I_{vf})+BC_{vwh}(I_{vfp}^s+I_{vfp}^w+I_{vf}^w))S_h}{N_h^2} \\ 
            -\mu_h-\theta-\dfrac{\alpha(BC_{vh}(I_{vfp}+I_{vf})+BC_{vwh}(I_{vfp}^s+I_{vfp}^w+I_{vf}^w))S_h}{N_h^2} \\ \theta \\
            \dfrac{BC_{hv}I_hS_{vf}^w}{N_h^2} \\ -\dfrac{BC_{hv}I_hS_{vf}^w}{N_h^2} \\
            \dfrac{BC_{hv}I_hS_{vf}}{N_h^2} \\ -\dfrac{BC_{hv}I_hS_{vf}}{N_h^2} \\ 0 \\ 0 \\
            \dfrac{BC_{hv}I_hS_{vfp}^w}{N_h^2} \\-\dfrac{BC_{hv}I_hS_{vfp}^w}{N_h^2} \\
            \dfrac{BC_{hv}I_hS_{vfp}^s}{N_h^2} \\ -\dfrac{BC_{hv}I_hS_{vfp}^s}{N_h^2} \\
            \dfrac{BC_{hv}I_hS_{vfp}}{N_h^2} \\ -\dfrac{BC_{hv}I_hS_{vfp}}{N_h^2} \\ 0 \\ 0
        \end{bmatrix},
    \end{split}
\end{equation*}}
\allowdisplaybreaks{\begin{equation*}
    \begin{split}
        J_4 = \tiny
        \begin{bmatrix}
            b_h+ \dfrac{(BC_{vh}(I_{vfp}+I_{vf})+BC_{vwh}(I_{vfp}^s+I_{vfp}^w+I_{vf}^w))S_h}{N_h^2} \\
            - \dfrac{(1-\alpha)(BC_{vh}(I_{vfp}+I_{vf})+BC_{vwh}(I_{vfp}^s+I_{vfp}^w+I_{vf}^w))S_h}{N_h^2} \\ 
            -\dfrac{\alpha(BC_{vh}(I_{vfp}+I_{vf})+BC_{vwh}(I_{vfp}^s+I_{vfp}^w+I_{vf}^w))S_h}{N_h^2} \\ -\mu_h \\
            \dfrac{BC_{hv}I_hS_{vf}^w}{N_h^2} \\ -\dfrac{BC_{hv}I_hS_{vf}^w}{N_h^2} \\
            \dfrac{BC_{hv}I_hS_{vf}}{N_h^2} \\ -\dfrac{BC_{hv}I_hS_{vf}}{N_h^2} \\ 0 \\ 0 \\
            \dfrac{BC_{hv}I_hS_{vfp}^w}{N_h^2} \\-\dfrac{BC_{hv}I_hS_{vfp}^w}{N_h^2} \\
            \dfrac{BC_{hv}I_hS_{vfp}^s}{N_h^2} \\ -\dfrac{BC_{hv}I_hS_{vfp}^s}{N_h^2} \\
            \dfrac{BC_{hv}I_hS_{vfp}}{N_h^2} \\ -\dfrac{BC_{hv}I_hS_{vfp}}{N_h^2} \\ 0 \\ 0
        \end{bmatrix},
    \end{split}
\end{equation*}}
\allowdisplaybreaks{\begin{equation*}
    \begin{split}
        J_5 = \tiny
        \begin{bmatrix}
            0 & -\dfrac{BC_{vwh}S_h}{N_h} & 0 & -\dfrac{BC_{vh}S_h}{N_h} & 0 & 0 & 0 \\
            0 & \dfrac{(1-\alpha)BC_{vwh}S_h}{N_h} & 0 & \dfrac{(1-\alpha)BC_{vh}S_h}{N_h} & 0 & 0 & 0 \\
            0 & \dfrac{\alpha BC_{vwh}S_h}{N_h} & 0 & \dfrac{\alpha BC_{vh}S_h}{N_h} & 0 & 0 & 0 \\
            0 & 0 & 0 & 0 & 0 & 0 &0\\
            -\mu_{fw}-\dfrac{BC_{hv}I_h}{N_h}-\sigma & 0 & 0 & 0 & 0 & 0&0 \\
            \dfrac{BC_{hv}I_h}{N_h} & -\mu_{fw}-\sigma & 0 & 0 & 0 & 0&0 \\
            0 & 0 & -\mu_{fu}-\dfrac{BC_{hv}I_h}{N_h}-\sigma & 0 & 0 & 0&0 \\
            0 & 0 & \dfrac{BC_{hv}I_h}{N_h} & -\mu_{fu}-\sigma & 0 & 0&0 \\
            0 & 0 & 0 & 0 & -\mu_{mw} & 0 &0 \\ 0 & 0 & 0 & 0 & 0 & -\mu_{mu}&0 \\
            \sigma & 0 &0&0&0&0&-\mu_{fw}-\dfrac{BC_{hv}I_h}{N_h} \\ \dfrac{BC_{hv}I_h}{N_h}&\sigma&0&0&0&0 \\
            0&0&m_w\sigma&0&0&0&0 \\ 0&0&0&m_w\sigma&0&0&0 \\
            0&0&m\sigma&0&0&0&0 \\ 0&0&0&m\sigma&0&0&0 \\
            0&0&0&0&0&0&\eta_w v_w \\ 0&0&0&0&0&0&\eta_w v
        \end{bmatrix},
    \end{split}
\end{equation*}}
and
\allowdisplaybreaks{\begin{equation*}
    \begin{split}
        J_6 = \tiny
        \begin{bmatrix}
            -\dfrac{BC_{vwh}S_h}{N_h} & 0 & -\dfrac{BC_{vwh}S_h}{N_h} & 0 & -\dfrac{BC_{vh}S_h}{N_h} & 0 & 0 \\
            \dfrac{(1-\alpha)BC_{vwh}S_h}{N_h} & 0 & \dfrac{(1-\alpha)BC_{vwh}S_h}{N_h} & 0 & \dfrac{(1-\alpha)BC_{vh}S_h}{N_h} &0&0 \\
             \dfrac{\alpha BC_{vwh}S_h}{N_h} & 0 & \dfrac{\alpha BC_{vwh}S_h}{N_h} & 0 & \dfrac{\alpha BC_{vh}S_h}{N_h} & 0 & 0 \\
             0 & 0 & 0 & 0 & 0 & 0 & 0 \\
             0 & 0 & 0 & 0 & 0 & b_f\psi & 0 \\
             0 & 0 & 0 & 0 & 0 & 0 & 0 \\
             0 & 0 & 0 & 0 & 0 & 0 & b_f\psi \\
             0 & 0 & 0 & 0 & 0 & 0 & 0 \\
             0 & 0 & 0 & 0 & 0 & b_m\psi & 0 \\
             0 & 0 & 0 & 0 & 0 & 0 & b_m\psi \\
             0 & 0 & 0 & 0 & 0 & 0 & 0 \\
             0 & 0 & 0 & 0 & 0 & 0 \\
             0 & -\mu_{fw}-\dfrac{BC_{hv}I_h}{N_h} & 0 & 0 & 0 & 0 & 0 \\
             0 & \dfrac{BC_{hv}I_h}{N_h} & -\mu_{fw} & 0 & 0 & 0 & 0 \\
             0 & 0 & 0 & -\mu_{fu}-\dfrac{BC_{hv}I_h}{N_h} & 0 & 0 & 0 \\
              0 & 0 & 0 & \dfrac{BC_{hv}I_h}{N_h} & -\mu_{fu} & 0 & 0 \\
             \eta_w v_w & 0 & 0 & 0 & 0 & -\mu_a-\psi & 0 \\
            \eta_w v & 0 & 0 & \eta & \eta & 0 & 0-\mu_a - \psi\\
        \end{bmatrix}.
    \end{split}
\end{equation*}}

\bibliographystyle{unsrt}
\bibliography{cas-refs}

\end{document}